\documentclass[11pt,mathpaper]{article}
\usepackage[margin=1in]{geometry}
\usepackage{graphicx}
\usepackage{mathtools}
\usepackage{comment}
\usepackage{calc}
\usepackage{amsmath, ulem}
\usepackage{amsxtra}
\usepackage{enumerate}
\usepackage[shortlabels]{enumitem}
\usepackage[thicklines]{cancel}
\usepackage{amssymb,wasysym}
\usepackage{amsfonts}
\usepackage{pifont}
\usepackage{amsthm}
\usepackage[sans]{dsfont}

\usepackage{mathrsfs}
\usepackage[numbers,sort&compress]{natbib}
\usepackage[colorlinks]{hyperref}

\graphicspath{{figday/}}

\newcommand{\sL}{\mathsf{L}}
\newcommand{\sE}{\mathsf{E}}
\newcommand{\sP}{\mathsf{P}}
\newcommand{\sH}{\mathsf{H}}
\newcommand{\BioS}{\mathcal K_{\mathsf{BS}}}
\newcommand{\bT}{\mathbf{T}}
\newcommand{{\bJ}}{\mathbf{J}}

\newcommand{\cVt}{\tilde{\mathcal{V}}}

\newcommand{\bP}{\mathbf{P}}
\newcommand{\bI}{\mathbf{Id}}

\newcommand{\bw}{\mathbf{w}}

\newcommand{\bZ}{\mathbf{Z}}

\newcommand{\bp}{\mathbf{p}}
\newcommand{\leps}{\bigg\rvert_{\varepsilon=0}}
\newcommand{\bomega}{\boldsymbol{\omega}}

\title{Variational Principle and Stochastic Lagrangian Formulation of Viscous Hydrodynamic Equations}
\author{Anping Pan\thanks{%
Department of Mathematics, Penn State University, 
University Park, PA 16802, USA. Email: abp5658@psu.edu}\quad Anna Mazzucato\thanks{%
Department of Mathematics, Penn State University, 
University Park, PA 16802, USA. Email: alm24@psu.edu}}

\begin{document}

\maketitle

\begin{abstract}
   In this manuscript, we extend the Lagrangian formulation of \cite{CI08} for Navier-Stokes Equation to a wider class of hydrodynamic models. Moreover, we prove that such Lagrangian formulation is naturally derived from a stochastic Hamilton-Pontryagin type variational principle. Generalized version of Kelvin circulation theorem in viscous fluids is discussed. We also derive self-contained local well-posedness results of some fluid equations based on Lagrangian-Eulerian formulation using fixed point argument. 
\end{abstract}

\tableofcontents

\section{Introduction}

 Variational formulation of inviscid hydrodynamic equations has a long history and can be traced back to the seminal work of Arnold \cite{Ar66} and Ebin-Marsden \cite{EM70}, where they established a formulation of incompressible Euler equation as geodesic in the infinite-dimensional Lie group of smooth volume-preserving diffeomorphisms, equipped with right-invariant $L^2$ Riemannian metric (kinetic energy). Since then generalizations to many other hydrodynamic PDEs have been obtained, we refer to \cite{AK98} \cite{KLM08} \cite{Wa16} for further discussions. Alternatively, from geometric mechanics viewpoint, in a unified abstract framework many inviscid hydrodynamic PDEs lie into the zoo of Euler-Poincar\'e equations (see \cite{MR99}) on  Lie algebras, and Lie-Poisson equations on dual Lie algebras as their Hamiltonian counterparts. This construction was extended to the class of fluid models with advected quantities using semi-direct product of Lie algebras, we refer to \cite{HMR98} for a detailed treatment on this topic. An in-depth discussion on geometric formulations of compressible fluids can be found in the recent survey \cite{KMM20}.

\paragraph*{Variational Principle for viscous Fluids} Given the variational formulation for inviscid fluids, it is natural to ask how to extend it to viscous hydrodynamics. Arnold’s framework is Newtonian (a second–order ODE) which does not directly accommodate viscosity in a geometric way. To bridge this inviscid–viscous gap, we consider the Lagrangian picture of linear advection-diffusion equation as a toy problem, where two alternative viewpoints arise. One is the Feynman-Kac formula, treating the advection-diffusion mechanism as averaging of tracers transporting along stochastic Lagrangian trajectories; Another is to view diffusion as effective gradient part of total drift, generated from Wasserstein gradient of entropy. Roughly speaking, the first perspective will lead to stochastic variational formulation and the second will lead to metriplectic formulation \cite{Mo86}. In present paper, we adopt the first perspective and extend the variational principle in \cite{Go05}, which-building on \cite{Os89}-introduces a magnetization variable to recast Navier-Stokes equation in dual Lie transport form, uncovering its nature of co-adjoint motion. We develop a general variational framework covering a range of viscous hydrodynamic PDEs and deriving a nonlinear Feynman-Kac representation that encompasses Constantin-Iyer's stochastic Lagrangian formulation of Navier-Stokes equation, which will be explained in detail in the next paragraph. 

We refer to \cite{CC07}\cite{Eyi10} for other attempts to generalize variational principles to Navier-Stokes setting. More recently, in \cite{Ho15} the author established variational principle for stochastic actions to derive stochastic hydrodynamic PDEs with multiplicative noise, while in \cite{ACC14} \cite{CCR23} the authors used mean drift of semi-martingales to obtain variational principles for deterministic Lagrangian. Our approach is, however, different from both.

\paragraph*{Lagrangian Formulation of Hydrodynamic Models}

 Lagrangian-Eulerian formulation of hydrodynamic equations was first derived by Constantin in \cite{Co01a} for Euler equation, and later by Constantin-Iyer in \cite{CI08} for Navier-Stokes equation, see also \cite{Co01b} \cite{CI10} \cite{FL18}\cite{Ya18} \cite{Zh10} \cite{Zh12}. Such formulation builds a coupled system in mixture of Eulerian velocity and Lagrangian flow and recovers the fluid velocity in Lagrangian form essentially using the dual Lie transport structure of fluid momentum equation. In case of Euler equation such system reads as:
\begin{equation}\label{e:LagEu}
\left\{
\begin{aligned}
&\dot X_t(x)=u(t,X_t(x))\\
&A_t(x)=X_t^{-1}(x)\\
&u(t,x)=\mathbf P[\nabla^*A_t(x)u_0(A_t(x))]
\end{aligned}
\right.
\end{equation}
where $A$ is the spatial inverse of flow map $X$, usually called the back-to-label map in the contexts of fluid dynamics. $\mathbf P$ is the Leray projection operator onto divergence free vector fields. Geometrically, Leray projector plays the role of cotangent-to-tangent bundle map, while evolution $\xi_t(x):=\nabla^*A_t(x)\xi_0(A_t(x))$ is the pull-back of co-vector under flow $X_t$.
 
In this framework,  the inviscid to viscous generalization is natural: Dissipative terms (Laplacian or fractional Laplacian) arise from averaging the scalar or vector moving along stochastic Lagrangian trajectories driven by Brownian or L\'evy noise. The resulting representation is therefore a Feynman-Kac type formula applied to advection-diffusion equations with McKean-type nonlinearity (see \cite{Mc66}), since the drift of stochastic flow depends on distribution of flow itself. For Navier-Stokes case in \cite{CI08}, such formulation is given by following system:
\begin{equation}\label{e:Lagrangian NSE}
\left\{
\begin{aligned}
&dX_t(x)=u(t,X_t(x))+\sqrt{2\nu}dW_t\\
&A_t(x)=X_t^{-1}(x)\\
&u(t,x)=\mathbb E_W\mathbf P[\nabla^*A_t(x)u_0(A_t(x))]
\end{aligned}
\right.
\end{equation}
where $W$ is the Wiener process on $\mathbb R^d$ and $\mathbb E_W$ is the expectation with respect to $W$.

Our main result in this work is to extend above stochastic Lagrangian formulations to a huge class of hydrodynamic models. We achieve this by realizing the geometric essence of \eqref{e:Lagrangian NSE} and finding a general Hamilton-Pontryagin type variational principle that naturally incorporates such stochastic Lagrangian systems as critical points of the variational principle. We therefore unifies a couple of viscous hydrodynamic models as critical points of a general stochastic action principle in section 3.2.

Notably, Lagrangian-Eulerian formulation of hydrodynamic PDEs has a natural advantage in well-posedness theory in the sense of prevention of loss of derivative (due to a novel regularized mechanism of Leray Projector, as we'll explain in section 2) and applicability of Picard iteration method. In \cite{CI08} the classical local well-posedness in H\"older spaces were recovered, we also refer to \cite{FL19}, \cite{LO23} for application of such formulation in nonlinear stochastic PDEs with transport noise. In section 4, we'll develop well-posedness result for hydrodynamic models using this approach, which is self-contained from our viewpoint. 

\paragraph{Remarks on Circulation Conservation}

One remarkable consequence of Constantin-Iyer's Lagrangian formulation is the generalization of Kelvin circulation theorem to viscous setting. The Kelvin circulation theorem is a consequence of particle relabeling symmetry in incompressible inviscid fluids, which says that circulation of smooth Euler flow along any material loop is conserved. For viscous case, from Eulerian viewpoint such conservation law is deteriorated due to the presence of viscous dissipation term. However, in \cite{CI08} the authors consider circulation of Navier-Stokes flow along loop transported by stochastic Lagrangian trajectories, then in average the conservation of circulation still holds true, uncovering the backward martingale property of circulation. For further study of this subject, we refer to \cite{DH18}

It turns out that conservation of circulation strongly relies on dual Lie transport structure of the equation, which, in our variational formulation, is compatible with the equation satisfied by conjugate momentum. Therefore we could establish Kelvin circulation theorem for various hydrodynamic models  in both inviscid and viscous case.  In this work we revisit this point, in section 3.5 we point out that circulation, together with other conservation laws due to particle relabeling symmetry, can be deduced from duality between dual Lie transport of momentum equation and Lie transport of Lagrangian line motions, which is an extension of the remark \cite{Dr22}.

\paragraph{Structure of the Paper}

This manuscript is organized as follows: In section 2, we introduce notations and conventions of our paper and introduce several technical lemmas for later use. In section 3, we introduce our variational principle for viscous fluid PDEs and derive a stochastic Lagrangian formulation from the variational principle, several examples are carefully studied. We also discuss how our formulation echoes Brenier's theory of generalized flows, and how the Lagrangian formulation leads to conservation laws due to particle-relabeling symmetry. In section 4, we prove local well-posedness result for hydrodynamic model using Picard Iteration in our Lagrangian framework.

\section*{Acknowledgements} The authors would like to thank Peter Constantin, Theodore Drivas, Gautum Iyer and Alexei Novikov for stimulating discussions and valuable suggestions. A. Mazzucato was partially supported by the US National Science Foundation grants DMS-1909103, DMS-2206453, DMS-2511023  and Simons Foundation Grant 1036502. A. Pan was partially supported the US National Science Foundation grants DMS-1909103 and DMS-2206453 (PI Mazzucato).

\section{Basic Set up}

 In this section, we collect some notations and review some background knowledge of the present paper. We also carry out several technical lemmas for later use.

\subsection{Some Notations and Preliminaries}

\quad

To avoid boundary technicalities, in this work we choose physical domain $D$ to be $\mathbb T^d$ or $\mathbb R^d$ for $d=2$ or $3$. However, most of the results hold independent of physical dimensions. Throughout the paper letter $u,v,w$ will be used for vector fields on $D$ and letter $X$ is used for (stochastic) Lagrangian flows. Subscript $t$ is always used to indicate time-dependence but not temporal derivative. However in some scenario to avoid messy notations we may drop the time subscript, which, hopefully, will not cause any confusion. For a given vector field $v$, we'll use the notation $\nabla^*v$ for the transpose of $\nabla v$: $(\nabla^*v)_{ij}=\partial_iv_j$.

We'll use boldface uppercase letters for matrices and some operators. Denote by $\mathbf I$ the identity matrix and ${\bJ}$ the standard symplectic matrix in dimension $2$, that is
$$ {\bJ}=\bigg[\begin{array}{cc}
    0 & -1 \\
     1& 0
\end{array}\bigg]$$

Let $v: D\to\mathbb R^d$ be a smooth vector field, the curl operator $\nabla\times$ is equicalent to take the antisymmetric part of gradient of $v$: $\nabla\times v=\nabla v-\nabla^*v$. Moreover we have
\begin{equation*}
\nabla v-\nabla^*v=\omega\bJ\quad\text{ when }d=2,\quad \nabla v-\nabla^*v\cong \omega\quad \text{ when }d=3   
\end{equation*}
Here when $d=2$ $\omega=\partial_1v_2-\partial_2v_1$ is a scalar, while in $d=3$ the above identification is in the sense of $\mathbf{A}(3)\cong \mathbb R^3$, here $\mathbf{A}(3)$ is the space of $3\times 3$ anti-symmetric matrices.

The celebrated Helmholtz-Hodge decomposition theorem tells that space of $L^2$ vector fields admits following representation of orthogonal sum:
\begin{equation*}
L^2(D;\mathbb R^d)=H_\sigma\oplus H_\nabla\oplus H_{\mathbf{h}}
\end{equation*}
where $H_\sigma$ and $H_\nabla$ are respectively $L^2$-closure of smooth divergence free vector fields and smooth gradient fields, $H_{\mathbf{h}}$ is the space of harmonic vector fields and is ignored in our torus case by fixing the zero-mean of vector fields.

We denote by $\mathbf P$ the Leray projection operator, which projects a $L^2$-vector field onto div-free (and mean-free when $D=\mathbb T^d$) vector fields. In closed form we have
\begin{equation}\mathbf P=\bI-\nabla \Delta_{D}^{-1}\nabla\cdot
\end{equation}
We denote by $\mathbf{\Pi}=\bI-\mathbf P$ the projection onto $H_\nabla$. $\mathbf P$ is a singular integral operator of order $0$ commuting with partial derivatives, thus bounded on H\"older and  Sobolev spaces.

The Biot-Savart operator $\BioS$ is an integral operator reads as follows:
\begin{equation*}
\BioS \omega=\nabla\times(-\Delta_D)^{-1}\omega
\end{equation*}
here when $d=2$, we abuse the notation for a moment to regard the operator $\nabla\times$ as $\nabla^\perp$. $\omega$ is a scalar/vector field when $d=2/d=3$ respectively. Moreover, we have the following:
\begin{equation*}
\BioS(\nabla\times v)=\bP v,\quad \nabla\times\BioS v=v.    
\end{equation*}

Let $\Phi:D\to D$ be a $C^1$-diffeomorphism, assume $\theta$, $v$ are respectively scalar/vector field. We denote by $\Phi_\sharp \theta$ ($\Phi_\sharp v$ resp.) the push-forward of $\theta$ ($v$ resp.) along $\Phi$, i.e.
\begin{equation}\label{e:pf}
 \Phi_\sharp \theta=\theta\circ \Phi^{-1},\quad \Phi_\sharp v=(\nabla \Phi v) \circ \Phi^{-1}.
\end{equation}

We'll use the notation $\mathscr L_u$ to denote Lie derivative along divergence-free vector field $u$. We recall that for scalar function $\theta$ we have $\mathscr L_u\theta=u\cdot\nabla\theta$, while for vector field $v$ we have $\mathscr L_uv=[u,v]=u\cdot\nabla v-v\cdot\nabla u$. We use the notation $\mathscr L_u^*$ for dual of $\mathscr L_u$ in vector transport:
\begin{equation*}
\mathscr L_u^* v=u\cdot\nabla v+\nabla^*u \cdot v    
\end{equation*}

We now introduce our notation for function spaces. Throughout our paper we'll work on Sobolev spaces, but parallel results in H\"older space setting would also hold $verbatim$. We use $s$ as Sobolev exponent and assume $s>1+d/2$. We define:
\begin{equation*}
  H_\sigma^s:=\overline{C_\sigma^\infty}^{H^s},\quad  C_\sigma^\infty:=\{v\in C^\infty(D;\mathbb R^d): \nabla\cdot v=0\} 
\end{equation*}
 With a little bit abuse of notation, we use $H^s(D)$ to denote the space of both $H^s$-scalar valued and vector valued function on $D$. 

We'll use $\lVert\cdot\rVert_{p}$ for $L^p$ norm, $\lVert\cdot\rVert_{H^s}$ for $H^s$-fractional Sobolev norm and $\lVert \cdot\rVert_{\dot H_s}$ for order $s$ homogeneous Sobolev norm, all in space variable. Define $\Lambda=(-\Delta)^{1/2}$, we have $\lVert f\rVert_{\dot H_s}=\lVert \Lambda^{s}f\rVert_{2}$.

One fact that will be useful in the following is the relation between $\lVert\nabla\times v\rVert_{H^s}$ and $\lVert\nabla v\rVert_{H^s}$. Indeed for smooth vector field $v$ we have
\begin{equation*}\label{Div-free Sobolev norm Identity}
\lVert\nabla v\rVert_{\dot H^s}^2=\lVert\nabla\cdot v\rVert_{\dot H^s}^2+\lVert\nabla\times v\rVert_{\dot H^s}^2.
\end{equation*}
Indeed, say in $d=3$ we have:
\begin{equation*}
-\Delta v=\nabla\times\nabla\times v-\nabla(\nabla\cdot v)\quad\Longrightarrow\quad \lVert\nabla v\rVert_2^2=\lVert\nabla\times v\rVert_2^2+\lVert\nabla\cdot v\rVert_2^2
\end{equation*}
Therefore, for divergence free vector fields, we have $\lVert\nabla v\rVert_{H^s}^2\sim\lVert\nabla\times v\rVert_{H^s}^2$.

We'll use the norm $\lVert\cdot\rVert_{C_tH_x^s}$ for $v\in C([0,T];H^s)$, defined as
\begin{equation*}
 \lVert v\rVert_{C_tH_x^s}=\sup_{0\le t\le T}\lVert v(t,\cdot)\rVert_{H^s}.   
\end{equation*}

We now introduce our notion for functional derivative. Given $\Psi: v\in H^s\to \Psi(v)\in \mathbb R$ and $u, \delta u\in H^s$, the functional derivative of $\Psi$, if exists, satisfies
\begin{equation}\label{Functional Derivative}
 \frac{d}{d\varepsilon} \bigg\rvert_{\varepsilon=0} \Psi(u+\varepsilon\delta u)=\big\langle \frac{\delta\Psi}{\delta v}(u),\delta u\big\rangle_2
\end{equation}

Apriori we have ${\delta\Psi}/{\delta v}\in \mathbf L(H^s;(H^s)^*)\cong \mathbf L_2(H^s\times H^s)$, here $\mathbf L(E;F)$ is the space of bounded linear maps between Banach spaces $E$ and $F$, and $\mathbf L_2$ denotes space of continuous bilinear form.

On the other hand, when $\Psi$ is a smooth function on $H_\sigma^s$,  Gateaux derivative ${\delta \Psi}/{\delta v}$ is defined  via \eqref{Functional Derivative}, hence apriori we have ${\delta\Psi}/{\delta u} \in \mathbf L(H_\sigma^s,(H_\sigma^s)^*)$, and $(H_\sigma^s)^*\cong H^{-s}/\sim$, where the equivalence relation $\sim$ is given by $f\sim g$ iff $f-g=\nabla \varphi$ for some distributional gradient $\nabla\varphi$, since for any $\delta u\in H_\sigma^s$, $\langle\delta u,\cdot\rangle$ as a linear function on $H^s$ vanishes on gradient fields. Here to fix the gauge, we choose the div-free representative, i.e. we set
\begin{equation*}
(H_\sigma^s)^*\cong\{v\in H^{-s}:\langle\nabla\varphi, v\rangle_2=0\quad\text{for all }\varphi\in H^{s+1}\}.   
\end{equation*}

We'll also need to consider variation of functional defined on Sobolev diffeomorphisms. We denote by $\mathscr D^s$ the space of invertible $H^s(D;D)$ maps with $H^s$ inverse.  For functional $F:\mathscr D^s\to \mathbb R$ we define its functional derivative $(\delta F/\delta\Phi)(X)$ as
\begin{equation*}
\langle \frac{\delta F}{\delta \Phi}(X),\delta X\circ X\rangle_2=\frac{d}{d\varepsilon}\bigg\rvert_{\varepsilon=0} F(X_\varepsilon),\quad \text{where }X_\varepsilon=(\bI+\varepsilon\delta X)\circ X.     
\end{equation*}
where $\delta X: D\to\mathbb R^d$ is an arbitrary smooth vector-valued map. Such choice is reasonable since for sufficiently small $\varepsilon$ we could guarantee $X_\varepsilon\in \mathscr D^s$ by inverse function theorem.

\subsection{Some Transport Estimates}

\quad 

We introduce some estimates in this subsection for later use. Throughout the rest of this paper, $C$ will denote a generic constant and may vary from line to line. 

Suggested by \cite{CI08}, Lagrangian formulation of hydrodynamic models usually involves a typical operator $\mathcal W$, called Weber operator, which we recall here:
\begin{equation}\label{Weber Operator}
\mathcal W:(u,\ell)\in H_\sigma^s\times H^s \to \mathbf P(\nabla^*\ell u)\in H_\sigma^s.
\end{equation}
The range of Weber operator is guaranteed to be in $H_\sigma^s$ since standard computation implies:
\begin{align*}
&\partial_i(\mathcal W(u,\ell))_j=\mathbf P(\partial_i (\partial_j \ell^k u_k))=\mathbf P(\partial_{ij}\ell^ku_k+\partial_j \ell^k\partial_i u_k)\\
&=\mathbf P (\partial_j \ell^k\partial_i u^k-\partial_i \ell^k\partial_j u_k)=\mathbf P(\underbrace{\partial_j(\partial_i\ell^ku_k)}_{=\nabla (\partial_i\ell\cdot u)}-\partial_i\ell^k\partial_j u_k+\partial_j \ell^k\partial_i u_k).
\end{align*}
By Sobolev embedding $H^s\hookrightarrow W^{1,\infty}$ for $s>1+d/2$ we have
\begin{align*}
 &\lVert \partial_i \mathcal W(u,\ell)\rVert_{H^{s-1}}\le C(\lVert \nabla \ell\rVert_{H^{s-1}}\lVert \nabla u\rVert_\infty+\lVert \nabla \ell\rVert_\infty\lVert \nabla u\rVert_{H^{s-1}})\le C\lVert \ell\rVert_{H^{s}}\lVert u\rVert_{H^{s}}.
\end{align*}
We are interested in Weber operator of special type $\mathcal W(u\circ \ell,\ell)$ and we investigate further on apriori estimates of this operator. First we have $L^2$-estimate:
\begin{equation*}
\lVert\mathcal W(u\circ \ell,\ell)\rVert_2\le \lVert\nabla \ell\rVert_\infty \lVert u\circ \ell\rVert_2.   
\end{equation*}
In case $\ell$ is Lebesgue measure-preserving we moreover have
\begin{equation}\label{WebL2}
\lVert\mathcal W(u\circ \ell,\ell)\rVert_2\le \lVert\nabla \ell\rVert_\infty\lVert u\rVert_2.    
\end{equation}
Moreover, for any given vector field $u$, use the fact that $\bP=\BioS \circ \nabla\times$ and \eqref{Div-free Sobolev norm Identity} one conclude: 
\begin{equation}\label{e:Curl of Weber}
\lVert\nabla \mathcal W(u\circ \ell,\ell)\rVert_{H^{s-1}}=\lVert\nabla\times\mathcal W(u\circ \ell,\ell)\rVert_{H^{s-1}}=\lVert\nabla\times(\nabla^*\ell u\circ \ell)\rVert_{H^{s-1}} .   
\end{equation}
For $d=2$, \eqref{e:Curl of Weber} reads
\begin{equation*}
\lVert\nabla \mathcal W(u\circ \ell,\ell)\rVert_{H^{s-1}}=\lVert \det(\nabla \ell)\omega\circ \ell\rVert_{H^{s-1}}. 
\end{equation*}
For $d=3$ \eqref{e:Curl of Weber} reads
\begin{equation*}
\lVert\nabla \mathcal W(u\circ \ell,\ell)\rVert_{H^{s-1}}=\lVert \mathsf{Cof}(\nabla\ell)^*\omega\circ \ell\rVert_{H^{s-1}}=\lVert\det(\nabla\ell)(\nabla\ell)^{-1}\omega\circ \ell\rVert_{H^{s-1}}.    
\end{equation*}
Therefore in case $\det(\nabla \ell)=1$, the following estimate holds for both $d=2$ and $3$:
\begin{equation}
 \lVert \mathcal W(u\circ \ell,\ell)\rVert_{H^{s}}=\lVert \mathcal W(u\circ \ell,\ell)\rVert_2+\lVert \ell_\sharp\omega\rVert_{H^{s-1}}\le \lVert\nabla \ell\rVert_\infty\lVert u\rVert_2+\lVert\ell_\sharp\omega\rVert_{H^{s-1}}.
\end{equation}
Now for $v_1,v_2\in H_\sigma^s$ and $\ell_1,\ell_2\in H^s$, consider 
$$\mathcal W(v_1,\ell_1)-\mathcal W(v_2,\ell_2)=\mathbf P((\nabla^*(\ell_1-\ell_2))v_1)+\mathbf P(\nabla^*\ell_2(v_1-v_2))$$
$$=-\mathbf P((\nabla^*v_1)(\ell_1-\ell_2))+\mathbf P(\nabla^*\ell_2(v_1-v_2)).$$
Therefore we conclude:
\begin{equation} \label{e:Difference Weber L2 Estimate}
    \lVert \mathcal W(v_1,\ell_1)-\mathcal W(v_2,\ell_2)\rVert_2\le C\lVert\nabla v_1\rVert_\infty\lVert \ell_1-\ell_2\rVert_2+C\lVert\nabla \ell_2\rVert_\infty\lVert v_1-v_2\rVert_2.
    \end{equation}
    and 
    \begin{equation*}
    \lVert \nabla[\mathcal W(v_1,\ell_1)-\mathcal W(v_2,\ell_2)]\rVert_{H^{s-1}}    
    \end{equation*}
\begin{equation} \label{e:Diffrence Weber Sobolev Estimate}
   \le C\lVert\nabla v_1\rVert_{H^{s-1}}\lVert \nabla(\ell_1-\ell_2)\rVert_{H^{s-1}}+C\lVert\nabla \ell\rVert_{H^{s-1}}\lVert \nabla(v_1-v_2)\rVert_{H^{s-1}}.
    \end{equation}
    
Now we recall some classical facts propagation of $H^s$-regularity in transport equation. The next lemma is the Lemma 4 of \cite{PR16}.
\paragraph{Lemma 2.1} Assume $s>1+d/2$, $U>0$ and $u\in C([0,T];H_\sigma^s)$ satisfies $\lVert \nabla u\rVert_{C_tH_x^{s-1}}\le U$, then there exists a constant $C$ such that the back-to-label map $A$ of $u$ satisfies
\begin{equation}\label{Back to Label Map Estimate}
\lVert \nabla A-\mathbf I\rVert_{C_tH_x^{s-1}}\le  \exp(CTU)-1
\end{equation}

The following quantitative version of propagation of high Sobolev regularity for transport equation will be useful. For its proof we refer to \cite{BCD11}.

\paragraph{Lemma 2.2} Let $u\in C([0,T];H_\sigma^{s})$ for some $s>1+d/2$ and $f\in H^r$, $0\le r\le s$. Assume there exists $U>0$ such that:
\begin{equation*}
 \lVert \nabla u\rVert_{C_tH_x^{s-1}}\le U   
\end{equation*}
Let $X\in C([0,T];H^{s})$ be the flow generated by $u$. Let $A$ be the back-to-label map, then there exists a constant $C$ such that for any $0\le t\le T$:
\begin{equation}
\lVert f\circ A_t\rVert_{H^r} \le [1+C\exp(CTU)]^{\lfloor r\rfloor+1}\lVert f\rVert_{H^r}.
\end{equation}

The next lemma is a variant of Lemma 3.6 in \cite{Zh10}.

\paragraph{Lemma 2.3} Let $X,\tilde X$ be flow map associated with $u,\tilde u\in C([0,T];H_\sigma^s)$, let $A,\tilde A$ be the corresponding back-to-label map, then we have
\begin{equation}\label{Gronwall X}
\lVert X_t-\tilde X_t\rVert_2\le \exp(C\sup_{0\le t\le T}\lVert \nabla u\rVert_{H^{s-1}})\int_0^t\lVert u_\tau-\tilde u_\tau\rVert_2d\tau .  
\end{equation}
\begin{equation}
\lVert A_t-\tilde A_t\rVert_2\le \exp(C\sup_{0\le t\le T}\lVert \nabla u\rVert_{H^{s-1}})\int_0^t\lVert u_\tau-\tilde u_\tau\rVert_2d\tau.   
\end{equation}

Proof: We notice that
$$X(t,x)=x+\int_0^t u(\tau,X(\tau,x))d\tau,\quad \tilde X(t,x)=x+\int_0^t \tilde u(\tau,\tilde X(\tau,x))d\tau$$
acting $A(t,x)$ and $\tilde A(t,x)$ on the right we observe:
$$A(t,x)=x-\int_0^t u(\tau,A(t-\tau,x))d\tau,\quad \tilde A(t,x)=x-\int_0^t \tilde u(\tau,\tilde A(t-\tau,x))d\tau$$
Therefore
$$\lVert A_t-\tilde A_t\rVert_2\le \int_0^t\lVert u(\tau,A(t-\tau,x))-\tilde u(\tau,\tilde A(t-\tau,x))\rVert_2d\tau$$
$$\le \int_0^t\lVert u(\tau,A(t-\tau,x))-\tilde u(\tau,A(t-\tau,x))\rVert_2+\lVert \tilde u(\tau,A(t-\tau,x))-\tilde u(\tau,\tilde A(t-\tau,x))\rVert_2d\tau$$
$$\le \int_0^t\lVert u_\tau-\tilde u_\tau\rVert_2 d\tau+\sup_{0\le \tau\le t}\lVert\nabla \tilde u_\tau\rVert_\infty \int_0^t\lVert A(t-\tau)-\tilde A(t-\tau)\rVert_2d\tau$$
Then Gronwall inequality implies
$$\sup_{0\le t\le T}\lVert A_t-\tilde A_t\rVert_2\le\exp(C\lVert \nabla u\rVert_{C_tH_x^{s-1}})\int_0^t \lVert u_\tau-\tilde u_\tau\rVert_2d\tau$$

The same argument implies estimate \eqref{Gronwall X} for $X$. \qed

\subsection{Lagrangian Formulas of Transport-Type Equations}

\quad

  A common geometric structure hidden in general incompressible hydrodynamic models is the following dual-Lie transport structure:
\begin{equation}\label{e:Dual Lie Transport}
 \partial_tw+u\cdot\nabla w+\nabla^*uw=f,\quad w(0,\cdot)=w_0.   
\end{equation}
Which is the dual of the following Lie-transport equation
\begin{equation}\label{e:Lie Transport}
 \partial_tv+u\cdot\nabla v-v\cdot\nabla u=0.   
\end{equation}

Notice \eqref{e:Lie Transport} preserves div-free vector fields, while \eqref{e:Dual Lie Transport} preserves gradient fields. manifests the co-adjoint geometric nature of many hydrodynamic equations. One note that given solution $w$ of \eqref{e:Dual Lie Transport}, the div-free part $\bw:=\bP w$ satisfies the following dual Lie transport equation:
\begin{equation}\label{e:Dual Lie Transport Projected}
\partial_t\bw+u\cdot\nabla \bw+\nabla^*u\bw+\nabla \bp=f,\quad \bw(0,\cdot)=\bP w_0    
\end{equation}
with some gradient field $\bp$. Indeed, let $w-\bw=\nabla q$, by \eqref{e:Dual Lie Transport} one has:
\begin{equation*}
\partial_t\bw+u\cdot\nabla \bw+\nabla^*u\bw+ \nabla(\underbrace{\partial_t q+u\cdot\nabla q}_{:=\bp})=f   
\end{equation*}
The following lemma summarizes Lagrangian representation of its solution and its behavior under taking curl operator, which can be viewed as a generalization of Cauchy vorticity formula. 

\paragraph{Lemma 2.4} Let $u\in C([0,T];H_\sigma^s)$,$f\in C([0,T];H^{s})$ and let $X$ be the flow of $u$, then the following holds:

\begin{itemize}
    \item {(i)} The unique solution $w_t$ of \eqref{e:Dual Lie Transport}
admits Lagrangian representation
\begin{equation}\label{e:Lagrangian Formula of Dual Lie Transport}
    w(t,x)=\nabla^*A_t(x)\bigg(w_0+\int_0^t \nabla^*X_\tau(x)f(\tau,X_\tau(x))d\tau\bigg)\circ A_t(x)
\end{equation}

\item (ii) In dimension $d=2$ ($d=3$, resp.) vorticity $\zeta:=\nabla\times w$ reduces to a scalar (vector, resp.) field. Moreover in Lagrangian form:
\begin{equation}\label{Vorticity Lagrangian Formula}
 \zeta(t,x)=\bigg[\nabla X_t\bigg(\nabla\times w_0+\int_0^t (A_\tau)_\sharp (\nabla\times f_\tau)d\tau\bigg)\bigg]\circ A_t(x)   
\end{equation}

\item (iii)  $\bw:=\BioS\zeta$ solves \eqref{e:Dual Lie Transport Projected}. Moreover we have $\zeta\in C_tH_x^{s-1}, \bw\in C_tH_x^s$ provided $w_0\in H^s$, $f\in C_tH_x^s$.

\end{itemize}

\quad

\textbf{Proof:} (i) Introduce Lagrangian vector field $\tilde w(t,x)=w(t,X_t(x))$, we observe 
\begin{align*}
 \frac{d}{dt}\tilde w(t,x)=(\partial_tw+u\cdot\nabla w)(t,X_t(x))=(f-\nabla^*uw)(t,X_t(x))   
\end{align*}

therefore $\tilde w(t,x)$ solves the following ODE
\begin{equation}\label{e:Nonhomogenous ODE}
\frac{d}{dt}\tilde w+\nabla^*u(X)\tilde w=f(X),\quad \tilde w_0=w_0   
\end{equation}

First note that the unique solution of
\begin{equation*}
\frac{d}{dt}\tilde w+\nabla^*u(X)\tilde w=0,\quad \tilde w_0=w_0   
\end{equation*}
is given by $\tilde w=(\nabla^*X_t)^{-1}w_0$. By Duhamel's principle, solution of nonhomogeneous ODE \eqref{e:Nonhomogenous ODE} is given by
\begin{align*}
  \tilde w(t,x)=(\nabla^*X_t)^{-1}(x)\bigg( w_0(x)+ \int_0^t\nabla^*X_\tau(x) f(\tau,X_\tau(x))d\tau\bigg) 
\end{align*}

Notice $(\nabla^*X_t)^{-1}=\nabla^*A_t(X_t)$, we then conclude $w(t,x)=\tilde w(t,A_t(x))$ is given by \eqref{e:Lagrangian Formula of Dual Lie Transport}. 

\quad

(ii) We now show that desired result holds for $d=2,3$. First let $d=2$, we notice 
\begin{equation*}
 \zeta(t,x)=\nabla\times [\nabla^*A_t(x)w_0\circ A_t(x)]  =\nabla\cdot{\bJ}^*[\nabla^*A_t(x)w_0\circ A_t(x)]  
\end{equation*}

\begin{equation*}
=\nabla\cdot [(\nabla A_t)^{-1}(x){\bJ}^* w_0\circ A_t(x)]  
=\nabla\cdot [(\nabla X_t{\bJ}^*w_0)\circ A_t(x)]  
\end{equation*}

\begin{equation*}
=\nabla\cdot [(X_t)_\sharp ({\bJ}^*w_0)]=( X_t)_\sharp(\nabla\cdot{\bJ}^*w_0)=(X_t)_\sharp(\nabla\times w_0).  
\end{equation*}
Which is our desired result.

To deal with the $d=3$ case, take any constant vector $\mathbf q\in \mathbb R^3$ we note that
\begin{equation*}
\zeta\times \mathbf q=(\nabla w-\nabla^*w)\mathbf q=\nabla^*A[(\nabla w_0-\nabla^*w_0)\circ A]\nabla A\mathbf q .   
\end{equation*}

Therefore we could rewrite
\begin{equation*}
 \zeta\times \mathbf q=\nabla^*A[(\nabla\times\omega_0)\circ A]\times (\nabla A\cdot\mathbf q ) . 
\end{equation*}

Now since for any $3\times 3$ matrix $\mathbf C$ and any vector $\mathbf a,\mathbf b\in\mathbb R^3$, we have
\begin{align*}
 (\mathbf C\mathbf a)\times (\mathbf C\mathbf b)=\det(\mathbf C)(\mathbf C^*)^{-1}(\mathbf a\times\mathbf b) . 
\end{align*}
Choosing $\mathbf C$ as $\nabla A$, $\mathbf a=(\nabla A)^{-1}\zeta(0, A)=(X_t)_\sharp\zeta(0,\cdot)$ and $\mathbf b=\mathbf q$ we have

\begin{equation*}
\zeta\times\mathbf q=\nabla^*A(\nabla^*A)^{-1}[(X_t)_\sharp\zeta_*(0,)]\times\mathbf q=[(X_t)_\sharp\zeta(0,\cdot)]\times\mathbf q    
\end{equation*}

By arbitrariness of $\mathbf q$ we conclude $\zeta(t,\cdot)=(X_t)_\sharp\zeta_0$ in the sense of vector field in dimension $3$, hence $\zeta$ solves \eqref{Vorticity Lagrangian Formula}.

\quad

(iii) To see that $\bw$ solves \eqref{e:Dual Lie Transport Projected}, recall that $\bP v=\BioS\nabla\times v$ for $v\in H^s$. Now choose $v$ as the solution $w$ of \eqref{e:Dual Lie Transport} such that $\nabla\times w_0=\zeta_0$, we immediately conclude that $\bw$ solves \eqref{e:Dual Lie Transport Projected}.

For the second part of the statement we'll show the following quantitative estimate: As in Lemma 2.1, take $U\ge \lVert\nabla u\rVert_{C_tH_x^{s-1}}$, we claim
\begin{equation} \label{e:Lie transport Estimate}
 \lVert\zeta\rVert_{C_tH_x^{s-1}}\le C[1+C\exp(CTU)]^{\lfloor s\rfloor+1}(\lVert  w_0\rVert_{H^{s}}+T\lVert \nabla\times f\rVert_{C_tH_x^{s-1}})   
\end{equation} 
 which naturally implies the desired result. We just deal with the case $d=3$, $d=2$ case $(X_t)_\sharp$ is pure rearrangement and thus desired result follows from Lemma 2.2. Note that 
\begin{equation*}
\lVert\zeta_t\rVert_{H^{s-1}}\le \underbrace{\lVert (\nabla X_t\nabla\times w_0)\circ A_t\rVert_{H^{s-1}}}_{:=(I)}+\underbrace{\bigg\lVert \int_0^t (X_{t-\tau})_\sharp(\nabla\times f_\tau)d\tau\bigg\rVert_{H^{s-1}}}_{:=(II)}
\end{equation*}
Now bound of (I) follows from lemma 2.2 and algebra structure of $H^{s-1}$:
\begin{equation*}
 (I)\le [1+C\exp(CTU)]^{\lfloor s\rfloor}\lVert \nabla X_t(\nabla\times w_0)\rVert_{H^{s-1}}\le C[1+C\exp(CTU)]^{\lfloor s\rfloor+1}\lVert  w_0\rVert_{H^{s}}.    
\end{equation*}
While (II) follows from the triangle inequality and above estimate:
\begin{equation*}
 (II)\le  \int_0^t\lVert(X_{t-\tau})_\sharp(\nabla\times f_\tau)\rVert_{H^{s-1}}d\tau\le CT[1+C\exp(CTU)]^{\lfloor s\rfloor+1}\lVert \nabla\times f\rVert_{C_tH_x^{s-1}}.  
\end{equation*}
Therefore, putting above together we conclude \eqref{e:Lie transport Estimate}. \qed

\quad

Above lemma can be restated as a commutation relation between curl operator, Lie transport and its dual. Indeed in case of $\nabla\cdot u=0$ and $u,v$ sufficiently smooth, we have
\begin{equation}\label{e:Commutation}
 \mathscr L_u\nabla\times v=\nabla\times \mathscr L_u^*v,\quad \BioS\mathscr L_u \nabla\times v=\bP \mathscr L_u^* v   
\end{equation}

We now introduce the Lagrangian description of curl of Lie-transported vector fields. (e.g. current of magnetic fields)
\paragraph{Lemma 2.5} Let $u\in C([0,T];H_\sigma^s)$ and $X$ be the flow of $u$. Let $v_t=(X_t)_\sharp v_0$ be a vector field on $D$, then $J=\nabla\times v$ admits the following Lagrangian representation: 
\begin{equation*}
J=(X_t)_\sharp\bigg(\nabla\times v_0+\int_0^t (A_\tau)_\sharp (\mathscr L_{v_\tau}\omega_\tau)+(A_\tau)_\sharp Q(u,v)d\tau\bigg)    
\end{equation*}
where 
\begin{equation*}
 \omega=\nabla\times u,\quad Q(u,v):=2\sum_{j=1}^d\nabla v_j\times \nabla u_j   
\end{equation*}
here the notation $\times$ is just the standard wedge product for $2$-d and $3$-d vectors.

\quad

\textbf{Proof. } Since $v_t$ solves \eqref{e:Lie Transport} with initial data $v(0,\cdot)=v_0$, we have
\begin{equation*}
 \partial_t v_t+u_t\cdot(\nabla v_t-\nabla^* v_t)-v_t\cdot(\nabla u_t-\nabla^* u_t)=v_t\cdot\nabla^* u_t-u_t\cdot\nabla^* v_t,\quad v(0,\cdot)=v 
\end{equation*}
while taking curl on above equation we notice
\begin{equation*}
\partial_t J+\mathscr L_uJ-\mathscr L_v\omega=2\bigg(\sum_{j=1}^d\nabla v_j\times \nabla u_j\bigg) 
\end{equation*}
therefore by \eqref{Vorticity Lagrangian Formula} the desired conclusion follows. \qed

\quad

The quadratic term $Q(u,v)$ is a mixture of first derivative of $u$ and $v$, therefore it's $a$ $priori$ $H^{s-1}$ provided $s>1+d/2$ due to the algebra structure of $H^{s-1}$. Hence by \eqref{e:Lie transport Estimate} we moreover have:
\begin{equation}
\lVert J\rVert_{C_tH_x^{s-1}}\le  C[1+C\exp(CTU)]^{\lfloor s\rfloor+1}(\lVert  v_0\rVert_{H^{s}}+T\lVert \nabla\times \mathscr L_v \omega\rVert_{C_tH_x^{s-1}}+CTU\lVert v\rVert_{H^s})      
\end{equation}

We now recall some basic facts about SDEs and stochastic flows, we'll focus on the case of Brownian noise, and our results also hold for L\'evy noise after a few modifications. We refer to  \cite{IW89} \cite{Ku97} \cite{Ku04} for extensive treatment on this topic.

Given drift $u\in C([0,T];H^s)$, the following SDE admits unique solution:
\begin{equation}\label{e:Brownian SDE}
dX_t^\nu=u(t,X_t^\nu)dt+\sqrt{2\nu} dW_t,\quad X(0,x)=x    
\end{equation}
which, in integral form, reads as
\begin{equation*}
X_t^\nu(x)=x+\int_0^t u(\tau, X_\tau^\nu(x))d\tau+\sqrt{2\nu}W_t,\quad 0\le t\le T    
\end{equation*}
Then $X_t^\nu$ forms a stochastic flow of $H^s$-diffeomorphisms. To see this one use the following trick: Let $X_t:=X_t^\nu-\sqrt{2\nu} W_t$, then \eqref{e:Brownian SDE} reads:
\begin{equation}\label{e:Brownian RDE}
 \frac{d}{dt}X_t=u(X_t+\sqrt{2\nu} W_t)=:u_\nu(t,X_t)   
\end{equation}
Since $u_\nu(t,x):=u(t,x+\sqrt{2\nu} W_t)$ is spatially $H^s$ and continuous in time, classical Cauchy-Lipschitz theory in ODE ensures that $X_t$ forms a family of diffeomorphisms, hence the same for $X_t^\nu$. Notice moreover $\nabla X_t^\nu=\nabla X_t$ satisfies the linearized ODE
\begin{equation*}
\frac{d}{dt}\nabla X_t=\nabla u_\nu(t,X_t)\nabla X_t    
\end{equation*}
and Liouville equation 
\begin{equation}\label{e:Liouville Equation}
 \det(\nabla X_t)=\exp\bigg(\int_0^t\nabla\cdot u_\nu(\tau,X_\tau)d\tau\bigg)   
\end{equation} 
which ensures that $X_t$, $X_t^\nu$ preserves Lebesgue measure almost surely provided $u$ is divergence-free.

We denote by $A_t$ the back-to-label map of $X_t$ and $A_t^\nu$ the back-to-label map of $X_t^\nu$. Define the random shift $\phi(t,x)=x-\sqrt{2\nu} W_t$, then immediately we have:
\begin{equation*}
 A_t^\nu(x)=A(t,x-\sqrt{2\nu}  W_t)=A_t\circ \phi_t (x)  
\end{equation*}

We now show the following main Lemma towards Lagrangian solution of dual Lie transport equation with dissipation.
\paragraph{Lemma 2.6} Assume $u\in C([0,T];H_\sigma^{s+1}(D;\mathbb R^d))$, $f\in C([0,T];H^{s+1})$ and $X$ solves SDE \eqref{e:Brownian SDE}. Let $A$, $X$, $X_t^\nu$, $A_t^\nu$ be defined as above. Then the following holds:

\begin{itemize}
    \item (I) For deterministic $\xi_0\in H^{s+1}$, consider the following ODE with random coefficient:
    \begin{equation*}
     \frac{d}{dt}\tilde{\xi}+\nabla^*u(X_t^\nu)\tilde\xi=f(X_t^\nu),\quad {\tilde\xi}_0=\xi_0   
    \end{equation*}
    then its unique pathwise solution is given by
    \begin{equation}\label{tilde xi Equation}
     \tilde\xi=(\nabla^*X_t^\nu)^{-1}\bigg(\xi_0+\int_0^t\nabla^*X_s^\nu f(s,X_s^\nu)ds\bigg)\in C_t^1H_x^s  
    \end{equation}

\item (II) Vector field
\begin{equation}\label{Lagrangian solution of PDE}
\bw_t=\mathbb E_W\bP\bigg[\nabla^*A_t^\nu\bigg(w_0+\int_0^t \nabla^*X_\tau^\nu f(\tau,X_\tau^\nu)d\tau\bigg)\circ A_t^\nu\bigg]    
\end{equation}
is the unique strong solution to following  diffusion equation:
\begin{equation}\label{Fractional Diffusion PDE}
 \partial_t\bw+u\cdot\nabla \bw+\nabla^* u\bw+\nabla \bp=f+\nu\Delta\bw,\quad \nabla\cdot \bw=0
\end{equation}
with initial data $\bw(0,\cdot)=\bP w_0\in H^{s+1}$. Moreover, we have $\bw\in C([0,T];H^{s+1})$.

\end{itemize}

\textbf{Proof.} (I) Recall the notation $\phi_t: x\to x-\sqrt{2\nu}  W_t$ and $\phi_t^{-1}:x\to x+\sqrt{2\nu}  W_t$, hence $X_t^\nu=\phi_t^{-1}\circ X_t$ and $A_t^\nu=A_t\circ \phi_t$. Consider random velocity field $u_\nu(t,x):=u(t,x+\sqrt{2\nu}  W_t)$. Unique solution of random matrix ODE:
    \begin{equation*}
    \dot{\mathbf Q}+\nabla^*u_\nu(X)\mathbf Q=0,\quad \mathbf Q_0=\bI \end{equation*}
is given by $\mathbf Q=(\nabla^*X)^{-1}=\nabla^*A\circ X$. Therefore the unique solution of following vector ODE with random coefficient
\begin{equation*}
\frac{d}{dt}\xi+\nabla^*u_\nu(X)\xi=0 \quad \xi(0,\cdot)=\xi_0   \end{equation*}
is given by $\xi_t:=(\nabla^*X)^{-1}\xi_0$. Here we could rewrite 
\begin{equation*}
\nabla^*u_\nu(X)=\nabla^*u(t,X+\sqrt{2\nu}  W_t)=\nabla^*u(X_t^\nu)   
\end{equation*}

Apply Duhamel's formula we obtain for nonhomogeneous ODE with force $f$:
\begin{equation*}
\frac{d}{dt}\xi+\nabla^*u_\nu(X)\xi=f(X_t^\nu) \quad \xi(0,\cdot)=\xi_0   \end{equation*}
the solution is given by
\begin{equation*}
\xi_t=(\nabla^*X_t)^{-1}\bigg(\xi_0+\int_0^t \nabla^* X_\tau f(X_\tau^\nu)d\tau\bigg)   
\end{equation*}

Note that here $\nabla X,\nabla A$ are in $C_t^1H_x^s$, hence $\xi\in C_t^1H_x^s$.

(II) Now we consider random vector field $\tilde w_t:=\xi_t\circ A_t$. Notice that
\begin{equation*}
 \tilde w_t=\nabla^* A_t\bigg(\xi_0+\int_0^t \nabla^* X_\tau f( X_\tau^\nu)d\tau\bigg)\circ  A_t   
\end{equation*}
which, by Lemma 2.4, is nothing but path-wise solution of the following random PDE:
\begin{equation}\label{e:RanPDE}
 \partial_t\tilde w+u_\nu\cdot\nabla \tilde w+\nabla^*u_\nu \tilde w=f\circ\phi_t^{-1}  
\end{equation}

 We now show that $\bw=\mathbb E_W\bP[\tilde w_t\circ \phi_t]$ is the unique strong  solution to \eqref{Fractional Diffusion PDE}. To this end, we first appreciate Lemma 2.4 to conclude that pathwisely:
 \begin{equation*}
   \tilde\bw:=\mathbf{P}\bigg[\nabla^* A_t\bigg(\xi_0+\int_0^t \nabla^* X_\tau f( X_\tau^\nu)d\tau\bigg)\circ  A_t\bigg]\in C([0,T];H^{s+1})   
 \end{equation*}
Now since $\tilde\bw$ has $C^2$ regularity, apply Ito-Wenztell formula to $\tilde\bw\circ \phi_t$ and average over the Brownian motion, the desired result follows from section 2.5, \cite{Iy06}. \qed

\paragraph{Remark 2.7} In this work we'll not bother ourselves to consider SDE with singular drift or driven by L\'evy noise. In case of singular drifts, severe technical issue arises when one apply the method of random characteristics. For studies concerning well-posedness of SDEs with drifts below Lipschitz regularity and associated parabolic PDEs in spirit of Feynman-Kac, we refer to \cite{KR05}\cite{Zh13}. For the case of L\'evy noise, similar result still holds but several technical modifications are necessary to made due to lack of time-regularity and momentum bound of L\'evy processes, we refer the readers to the thesis \cite{Pa26} for further discussions.

\section{Stochastic Variational Principle for Viscous Fluids}

We now develop a Hamilton-Pontryagin type variational principle (cf. \cite{BCHM00} \cite{HSS09}) for viscous hydrodynamic PDEs. The idea, as we explained above, is to perturb the deterministic Lagrangian by noise at the level of Lagrangian trajectory.  We'll derive various viscous hydrodynamic PDEs as critical points of our stochastic functional in this setting. Moreover, a generalization of Constantin-Iyer's formulation for a large family of viscous fluid PDEs is established. We'll also take advantage of this Lagrangian viewpoint to discuss conservation laws from particle relabeling symmetry. Finally, we'll discuss the connection between our variational principle and Brenier's generalized flow, which also relates to the groupoid description in \cite{IK24}.

\subsection{Geometry of Lagrangian-Eulerian Formulation of Euler's equation}

We now briefly explain the geometric insight of formulation \eqref{e:LagEu} in \cite{Co01a}, which will motivate our variational formulation of viscous fluids.

We recall the following form of Euler's equation in magnetization variable (cf. \cite{Go05}\cite{Os89}):
\begin{equation*}
\partial_t\xi+u\cdot\nabla\xi+\nabla^*u\xi+\nabla q=0,\quad u=\bP\xi    
\end{equation*}
Geometrically, the dual Lie transport equation lives in the dual-Lie algebra $\mathfrak X_\sigma^*$, the dual space of Lie algebra $\mathfrak X_\sigma$ of smooth div-free vector fields. The Leray projector $\bP$ formally plays the role of inertia operator (bundle map), sending momentum to velocity. Via the cotangent lift momentum map (see \cite{MR99}):
\begin{equation*}
J: (\mathbf{m},X)\to \mathbf{m}\circ X^{-1}    
\end{equation*}
we may lift back the $\xi$-dynamics to the cotangent bundle of measure-preserving diffeomorphism group. The lifted dynamics is the Hamiltonian ODE associated with the following Hamiltonian:
\begin{equation*}
H(u,X,\mathbf{m})=-\frac{1}{2}\lVert u\rVert_2^2+\langle \mathbf{m}, u(X)\rangle_2    
\end{equation*}
i.e., $(u,X,\mathbf{m})$ solves the following Hamiltonian ODE with dynamical constraint:
\begin{equation}\label{Hamiltonian ODE of HP}
\left\{
\begin{aligned}
&\dot X=\frac{\partial H}{\partial\mathbf{m}}\\
&\dot{\mathbf m}=-\frac{\partial H}{\partial X}\\
& \frac{\partial H}{\partial u}=0
\end{aligned}
\right.
\end{equation}
which corresponds to critical point of the following Hamilton-Pontryagin Lagrangian action:
\begin{equation*}
L(u,X,\mathbf{m})=\int_0^T\frac{1}{2}\lVert u\rVert_2^2-\langle\mathbf{m},\dot X-u(X)\rangle_2 dt    
\end{equation*}

Simple calculation verifies that the Euler-Lagrange equation of $L$ is equivalent to \eqref{Hamiltonian ODE of HP}. More precisely, we can show that \eqref{Hamiltonian ODE of HP} reads
\begin{equation*}
\left\{
\begin{aligned}
&\dot X=u(X)\\
&\dot{\mathbf m}+\nabla^*u(X)\mathbf{m}=0\\
& u=\bP[\mathbf m\circ X^{-1}]
\end{aligned}
\right.
\end{equation*}
Now, we can recover Constantin's Lagrangian-Eulerian formulation  \eqref{e:LagEu} of Euler using the above system, once we use the explicit Lagrangian formulae $\mathbf{m}=(\nabla^*X_t)^{-1}\mathbf{m_0}$ to represent the momentum dynamics $\mathbf{m}$. Indeed, assume $\mathbf{m}_0=u_0+\nabla q_0$, we have:
\begin{equation*}
u=\bP[\mathbf m\circ X^{-1}]=\bP[((\nabla^*X_t)^{-1}\mathbf{m_0})\circ X^{-1}]=\bP[\nabla^*A_tu_0\circ A_t]+\bP[{\nabla^*A_t\nabla q_0\circ A_t}]=\bP[\nabla^*A_tu_0\circ A_t].   
\end{equation*}
The last term vanishes, since we notice
\begin{equation*}
 \nabla^*A_t\nabla q_0\circ A_t=\nabla(q_0\circ A_t)   
\end{equation*}
Hence, we claim that $u=\bP[\mathbf m\circ X^{-1}]$ recovers the Weber's formula for incompressible Euler in \eqref{e:LagEu}.

We now summarize that the Lagrangian-Eulerian formulation, in general, is essentially the Lagrangian solution formulae of some infinite-dimensional Hamiltonian ODE on the cotangent bundle, therefore has a variational origin.

Recall that the Constantin-Iyer stochastic Lagrangian formulation is a stochastic perturbation of \eqref{e:LagEu} at Lagrangian level. Since in the inviscid case, \eqref{e:LagEu} is equivalent to critical point equation of a  Hamilton-Pontryagin type Lagrangian action, convincingly we should hear the echo of the stochasric Lagrangian formulation from the variational side, which is the stochastic variational principle that we'll discuss in the next subsection.

\subsection{Stochastic Hamilton-Pontryagin Principle}
Let $W_t$ be a $d$-dimensional Brownian motion. We introduce the following spaces: Let
\begin{equation*}
 \tilde{\mathscr D}:=L^\infty(\Omega,\mathbb P;C^1([0,T];\mathscr D^{s+1})),\cVt:=L^2(\Omega,\mathbb P;C^1([0,T];H^{s}))   
\end{equation*}
 Let $e_\omega$ be the evaluation map at $\omega\in\Omega$, therefore for any $\xi\in \cVt$ and $\omega\in\Omega$, we have $e_\omega(\xi)=\xi(\omega,\cdot)\in C^1([0,T];H^{s})$. 

Following previous notations, we denote by $X_t^\nu:=X_t+\sqrt{2\nu} W_t$. To avoid ambiguity, we use bold symbol $\bomega$ for vorticity and standard $\omega$ for randomness. 

We now state our main theorem:
\paragraph{Theorem 3.1} Let $\sE$ be a smooth function on $H_\sigma^{s+1}$ with specific form $\sE(v)=\frac{1}{2}\langle \mathbf{T}v,v\rangle_2$, where $\mathbf T: H^{s+1}\to L^2$ is a symmetric positive definite linear operator (a.k.a inertia operator) commuting with partial derivatives (e.g. Fourier multipliers), and let $\sP:\tilde{\mathscr D}\to \mathbb R$ satisfies the following condition:
\begin{equation*}
\frac{\delta\sP}{\delta\bZ}(X)\in L^2(\Omega\times [0,T];H_x^{s}) \quad\text{ for all } X\in \tilde{\mathscr D}
\end{equation*}

Let  $\Psi:(u,X,\tilde\xi)\in C_tH_\sigma^{s+1}\times \tilde{\mathscr D}\times \cVt \to \mathbb R$ be a functional defined as follows:
\begin{equation}\label{e:SFunc}
\Psi(u,X,\tilde\xi)=\int_0^T \sE(u_t)-\sP(X_t)+\mathbb E_W\langle\tilde\xi_t,\dot X-u(X_t^\nu)\rangle_2dt    
\end{equation}
Then the following holds:

\begin{itemize}
\item[(I)] (Euler-Lagrange Equation) The Euler Lagrange equation of \eqref{e:SFunc} reads as:
\begin{equation}\label{e:Stochastic Critical Point}
\left\{
\begin{aligned}
&\dot X=u(X+\sqrt{2\nu} W)=u(t,X_t^\nu)\\
&A_t^\nu=(X_t^\nu)^{-1}\\
&\dot{\tilde \xi}+\nabla^*u(t,X_t^\nu)\tilde\xi=-\frac{\delta \sP}{\delta\bZ}(X_t^\nu)\\
&\mathbf{T}u_t=\mathbb E_W\mathbf{P}(\tilde\xi_t\circ A_t^\nu)
\end{aligned}
\right.
\end{equation}
Hence, the critical process $X$ satisfies the condition that $X_t^\nu$ solves \eqref{e:Brownian SDE}.

\item[(II)] (Hamiltonian Formulation) Introduce the following random Hamiltonian:
\begin{equation}\label{e:RanHam}
\sH(\tilde\xi, X,u):=\mathsf E(u)+\mathsf P(X_t)-\langle \tilde\xi, u(X_t^\nu)\rangle_2   
\end{equation}
where $(\tilde\xi,X)$ is viewed as momentum-position variable and $u$ is treated as a dynamical constraint,  then \eqref{e:Stochastic Critical Point} is equivalent to the following random Hamiltonian ODE:
\begin{equation}\label{e:HamODE}
\left\{
\begin{aligned}
&\dot X=-{\delta \sH}/{\delta{\tilde\xi}}\\
&{\dot{\tilde\xi}}={\delta \sH}/{\delta X}\\
&\int_\Omega \frac{\delta \sH}{\delta u} \mathbb P(d\omega)=0
\end{aligned}
\right.
\end{equation}

\item[(III)] (Generalized Constantin-Iyer Formulation) Under initial condition $u(0,\cdot)=u_0$, $\tilde\xi_0=\mathbf Tu_0$ and $X_0=\bI$, with $u_0\in H^{s+1}$. For solution $(u,X,\tilde\xi)\in C_tH_\sigma^{s+1}\times \tilde{\mathscr D}\times \cVt$ of Hamiltonian ODE \eqref{e:HamODE}, we have that $(u,X)$ solves the following implicit Lagrangian-Eulerian system
\begin{equation}\label{e:Stochastic Lagrangian System}
\left\{
\begin{aligned}
&dX_t^\nu(t,x)=u(t,X_t^\nu(x))dt+\sqrt{2\nu} dW_t,\quad  A_t^\nu(x)=(X_t^\nu)^{-1}(x)\\
&u(t,x)={\mathbf T}^{-1}\mathbf P\mathbb E_W\bigg[\nabla^*A_t^\nu(x)\bigg(\mathbf Tu_0+\int_0^t \nabla^*X_\tau^\nu \frac{\delta \mathsf P}{\delta \mathbf Z}(X_\tau^\nu) d\tau\bigg)\circ A_t^\nu(x)\bigg].
\end{aligned}
\right.
\end{equation}

Conversely, given $u_0\in H_\sigma^{s+1}$, assume the  fixed point $(u,A^\nu,X^\nu)$ of \eqref{e:Stochastic Lagrangian System} satisfies $u\in C_tH_\sigma^{s+1}$, then $(u,A^\nu,X^\nu)$ produces solution of $(u,X,\tilde\xi)$ of \eqref{e:HamODE} by setting $X_t=X_t^\nu-\sqrt{2\nu}W_t$, and setting $\tilde\xi$ by \eqref{tilde xi Equation} with $f=\delta\sP/\delta\bZ$ as follows: 
\begin{equation}\label{tilde xi Equation II} 
\tilde\xi=(\nabla^*X_t^\nu)^{-1}\bigg(\xi_0+\int_0^t\nabla^*X_\tau^\nu \frac{\delta\sP}{\delta \bZ}(\tau,X_\tau^\nu)d\tau\bigg).    
\end{equation}

\item [(IV)] (Eulerian Formulation) Assume that $(u,X^\nu,A^\nu)$ is a solution of \eqref{e:Stochastic Lagrangian System} with initial data $u_0$ such that $u\in C_t^1H_\sigma^s$, let $\tilde\xi$ be defined as 
\eqref{tilde xi Equation II} and moreover define
\begin{equation*}
 \xi=\mathbb E_W\bP[\tilde\xi\circ A_t^\nu].   
\end{equation*}
Then $(u,\xi)$ solves the following Cauchy problem of PDE system:
\begin{equation}\label{e:Eulerian Momentum}
\left\{
\begin{aligned}
&\partial_t\xi+u\cdot \nabla\xi+\nabla^*u\xi+\mathbb E_W\big[\frac{\delta \sP}{\delta\bZ}(X_t)\circ A_t^\nu\big]+\nabla q=\nu\Delta\xi\\
&u=\bT^{-1}\xi\\
&\nabla\cdot\xi=0\\
&\xi(0,\cdot)=\bT u_0,\quad u(0,\cdot)=u_0
\end{aligned}
\right.
\end{equation}

\item[(V)] (Vorticity Formulation)  Assume $(u,\xi)$ is a solution of \eqref{e:Eulerian Momentum}. Then for $\bomega:=\nabla\times \xi$, $(u,\bomega)$ solve the following active scalar/vector transport equation:

\begin{equation}\label{e:Vorticity Formulation}
\left\{
\begin{aligned}
&\partial_t\bomega+\mathscr L_u\bomega=\nu\Delta\bomega-\nabla\times\mathbb E_W(\frac{\delta \sP}{\delta\bZ}(X_t^\nu)\circ A_t^\nu)\\
&\mathbf{T}u={\BioS}\bomega
\end{aligned}
\right.
\end{equation}
Conversely, given solution $(u,\bomega)$ of \eqref{e:Vorticity Formulation},  we have $\xi:={\BioS}\bomega$ is the solution of \eqref{e:Eulerian Momentum}.

\end{itemize}

We remark that for some natural choice of combination of $\bT$ and $\sP$, the fixed point $u$ of \eqref{e:Stochastic Lagrangian System} satisfies $u\in C_tH_\sigma^{s+1}$ provided $u_0\in H_\sigma^{s+1}$, we delay the proof of these examples to section 4. 

\quad

\textbf{Proof of Theorem 3.1}:  (I)  We start with computing the first variation of functional $\Psi$ with respect to variables $u$, $X$ and $\tilde \xi$. Consider compact perturbation $\delta u\in H_\sigma^{s}$, $\delta X\in L^2(\Omega;C_c^1([0,T];H_x^s))$ and $\delta\tilde\xi\in \cVt$ of $\Psi$ respectively. Under above choice of $\tilde{\mathscr D}$, $\cVt$, all the pairings and integrals below are well-defined.  Now we compute corresponding variational equations as follows:
\begin{itemize}
\item [(i)] The variation in variable $\tilde\xi$ simply reads as
\begin{equation}
\dot X=u(X+\sqrt{2\nu} W_t)=u(X_t^\nu) 
\end{equation}

which implies that $X_t^\nu$ is a stochastic flow with drift $u$. Therefore $X_t^\nu$ is a.s. invertible and Lebesgue measure preserving, since we have by \eqref{e:Liouville Equation}
\begin{equation*}
\det(\nabla X_t^\nu)=\det(\nabla X_t)=\exp\bigg(\int_0^t\nabla\cdot u(\tau,X_\tau+\sqrt{2\nu} W_\tau)d\tau\bigg)=1    
\end{equation*}
guarantees that $\det(\nabla A_t^\nu)=1/\det(\nabla X_t^\nu\circ A_t^\nu)=1$, $\mathbb P$-a.s.

\item [(ii)] The variation in variable $u$ is given by

\begin{equation*}
\frac{d}{d\varepsilon}\bigg\rvert_{\varepsilon=0}\Psi(u+\varepsilon\delta u,X,\tilde\xi)=\int_0^T \big\langle \frac{\delta \sE}{\delta v}(u_t),\delta u_t\big\rangle_2-\mathbb E_W\langle \tilde\xi_t, \delta u(X_t^\nu)\rangle_2 dt    
\end{equation*}
Using the change of variable formula, vanishing of RHS for all $\delta u\in C^1([0,T];H_\sigma)$ implies existence of scalar function $q$ such that 
\begin{equation*}
 \bT u=\frac{\delta \sE}{\delta v}(u_t)=\mathbb E_W(\det(\nabla A_t^\nu)\tilde\xi\circ A_t^\nu)+\nabla q  
\end{equation*}
Moreover since $\mathbf{T}$ commute with partial derivatives, we have $\nabla\cdot \mathbf{T}u=0$. Therefore we conclude as expected
\begin{equation*}
u=\bT^{-1}\mathbf P\mathbb E_W[\det(\nabla A_t^\nu)\tilde\xi\circ A_t^\nu]    
\end{equation*}
Combining with  (i) we conclude critical $u$ satisfies
\begin{equation}\label{Critical u}
u=\bT^{-1}\mathbf P\mathbb E_W(\tilde\xi_t\circ A_t^\nu)    
\end{equation}
\item [(iii)] We now compute variation in $\delta X$. 
\begin{align*}
 &\frac{d}{d\varepsilon}\bigg\rvert_{\varepsilon=0}\Psi(u,X+\varepsilon\delta X\circ X,\tilde\xi)\\
 &=\int_\Omega\int_0^T -\langle \frac{\delta \sP}{\delta \mathbf Z}(X_t),\delta X_t\circ X_t\rangle_2 dt\mathbb P(d\omega)+\mathbb E_W\langle \tilde\xi_t,\frac{d}{dt}({\delta X_t\circ X_t})-\nabla u(X_t^\nu)\delta X_t\circ X_t\rangle  dt\\
 &=\int_0^T -\langle \frac{\delta \sP}{\delta \mathbf Z}(X_t),\delta X_t\circ X_t\rangle_2-\mathbb E_W\langle \partial_t\tilde\xi_t+\nabla^*u(X_t^\nu)\xi,\delta X_t\circ X_t\rangle_2 dt\\
 &=-\int_\Omega\int_0^T\int_{D}(\delta X_t\circ X_t)\cdot\bigg(\underbrace{\dot{\tilde \xi}+\nabla^*u(t,X_t^\nu)\tilde\xi+\frac{\delta \sP}{\delta\mathbf Z}(X_t^\nu)}_{(*)}\bigg)dxdt\mathbb P(d\omega)
\end{align*}

where from second to third line we integrated by part in time and use compactly-support in time of variation $\delta X$. Now by assumption on $u,X,\tilde\xi$ and $\sP$, apriori we have 
\begin{equation*}
\dot{\tilde\xi},\frac{\delta\sP}{\delta\bZ}(X_t)\in L^2(\Omega\times [0,T];H_x^{s}),\quad \nabla^*u_\nu(X)\in L^\infty(\Omega\times [0,T];H_x^{s})
\end{equation*}
Thus by arbitrariness of $\delta X$, it's necessary to guarantee $(*)=0$ to make the $X$-variation vanishes. Hence, we obtain \eqref{e:Stochastic Critical Point} by putting above together.

\end{itemize}

(II) For random Hamiltonian \eqref{e:RanHam}, similar computation of variation yields: for each fixed path $\omega\in\Omega$ we have:
\begin{equation*}
\frac{\delta \sH}{\delta\tilde\xi}(\omega)=u(X_t^\nu(\omega)),\quad \frac{\delta \sH}{\delta X}(\omega)=-\nabla^*u(X_t^\nu(\omega))\tilde\xi_t(\omega)+\frac{\delta \sP}{\delta\bZ}(X_t(\omega)).
\end{equation*}
Meanwhile for variational derivative in $u$ it reads:
\begin{equation*}
\frac{\delta \mathsf H}{\delta u}(\omega,t,x)= \bP [(\bT u)(t,x)-\tilde\xi_t(\omega,A_t^\nu(\omega,x))].
\end{equation*}
Therefore \eqref{e:HamODE} immediately follows.

\quad

(III) Assume now $(u,X,\xi)$ solves the Hamiltonian equation \eqref{e:RanHam}. By the first equation we immediately conclude $X_t^\nu$ satisfies the SDE. To recover $u$, simply notice that $u$ is recovered from $\tilde\xi$ via \eqref{Critical u} and apply Lemma 2.6 to $\mathbb E\bP (\tilde\xi\circ A_t^\nu)$. 

For the converse direction, appealing to the assumption $u\in C_tH_\sigma^{s+1}$, we are able to apply lemma 2.6, which demonstrates that
$\tilde\xi$ defined as \eqref{tilde xi Equation II}, solves
\begin{equation*}
 \dot{\tilde\xi}+\nabla^*u(X_t^\nu)\tilde\xi=\frac{\delta\sP}{\delta\bZ}(X_t^\nu).   
\end{equation*}
Together with \eqref{Critical u} the desired result follows.

\quad

(IV) We now show \eqref{e:Eulerian Momentum}. Notice that the random ODE for $\tilde \xi$ can be rewritten as
\begin{equation*}
 \dot{\tilde\xi}+\nabla^*u_\nu(X_t){\tilde\xi}=-\frac{\delta \sP}{\delta\mathbf Z}(X).
\end{equation*}
Now write $\xi=\mathbb E_W\tilde\xi\circ A_t^\nu=\mathbb E_W(\xi\circ A_t\circ \phi_t^\nu)$, \eqref{e:Eulerian Momentum} follows from Lemma 2.6 (ii).

\quad

(V) Let $\omega=\nabla\times\xi=\nabla\times ({\delta \sE}/{\delta v})(u_t)$,  since partial derivatives commute with $\Delta$, the proof follows from Lemma 2.4 where the follow commutation relation was shown:
\begin{equation*}
\nabla\times(u\cdot\nabla \xi+\nabla^*u\xi)=\mathscr L_u(\nabla\times\xi).    
\end{equation*}
Conversely, assume $(u,\omega)$ solves \eqref{e:Vorticity Formulation}, then due to commutativity of Biot-Savart operator and $\Delta$, together with 
lemma 2.4 and lemma 2.6,  the proof is completed.

\qed

\paragraph{Remark 3.2} It's straightforward to see that dropping randomness in our setting, we will obtain Hamilton-Pontryagin principle and Lagrangian-Eulerian formulation of many inviscid hydrodynamic models. More precisely, Introduce augmented Lagrangian $\Phi$ on $C([0,T];H_\sigma^s)\times C^1([0,T];\mathscr D^s)\times C([0,T];H^{s-1})$:
\begin{equation}\label{e:Inviscid Augmented Lagrangian}
\Phi(u,X,\xi)=\int_0^T F(u_t,X_t)-\langle\xi_t,u_t-\dot X_t\circ X_t^{-1}\rangle_2 dt.
\end{equation}
Then following the same line as in theorem 3.1, we get critical point equation of $\Psi$ and corresponding Eulerian/vorticity/Lagrangian Eulerian formulations, which includes various inviscid fluid equations by choosing appropriate $F$. Here more general structure of $F$ is allowed so we can include examples like inhomogeneous incompressible Euler (See \cite{Da06}) or $p$-Euler equation (See \cite{LL18}). For more details we refer to the ongoing thesis \cite{Pa26}.

\subsection{Examples}

We now list several examples of viscous PDEs that could be derived from above variational principle, and their associated Lagrangian formulations.

\paragraph{Example 3.3} Consider the $2$-d example 
\begin{equation*}
 \sE(u)=\frac{1}{2}\langle (-\Delta)^{-\alpha}u, u\rangle_2,\quad \sP=0  
\end{equation*}
where $0\le \alpha\le 1$. Then the corresponding critical equation in vorticity form \eqref{e:Vorticity equation} is:
\begin{equation*}
\partial_t\omega+u\cdot\nabla\omega=\nu\Delta\omega,\quad u=\nabla^\perp(-\Delta)^{-1+\alpha}\omega
\end{equation*}

which is the generalized SQG equation with dissipation. When $\alpha=0$ we obtain Navier-Stokes equation as in \cite{CI08}, while for $\alpha=1$, we find that $\omega$ solves the heat equation. The Lagrangian formulation reads:
\begin{equation*}
\left\{
\begin{aligned}
&dX_t^\nu(t,x)=u(t,X_t^\nu(x))dt+\sqrt{2\nu} dW_t,\quad  A_t^\nu(x)=(X_t^\nu)^{-1}(x)\\
&u(t,x)=(-\Delta)^{\alpha}\mathbf P\mathbb E_W[\nabla^*A_t^\nu(x)((-\Delta)^{-\alpha}u_0\circ A_t^\nu(x))]
\end{aligned}
\right.
\end{equation*}

We can derive the viscous Camassa-Holm equation and other generalized viscous fluid models (e.g. see \cite{CIW08}, \cite{CCCGW12}) by simply modifying the Fourier multiplier in constitutive law. We can also replace Brownian noise with L\'evy noise to derive Generalized SQG with fractional dissipation.

\quad

We now move on to discuss Lagrangian with potential, which leads to fluid-tracer coupling models with unit "Prandtl" number. We refer to \cite{Eyi09} for an equivalent but slightly different stochastic Lagrangian formulation of MHD system with the same viscosity and resistivity.
\paragraph{Example 3.4} Consider the energy functional
\begin{equation*}
 \sE=\frac{1}{2}\lVert u\rVert_2^2,\quad \sP(X_t)=\underbrace{-\mathbb E_W\langle \theta_0,\varphi(X_t^\nu)\rangle_2}_{:=\sP_1(X)}+\underbrace{\frac{1}{2}\lVert\mathbb E[(X_t^\nu)_\sharp B_0] \rVert_2^2}_{\sP_2(X)}
\end{equation*}
here $\varphi$ is a given smooth potential function. Then we have the variation
\begin{equation*}
\frac{d}{d\varepsilon}\leps\sP_1(X+\varepsilon\delta X\circ X)=-\mathbb E_W\langle \nabla\varphi(X_t^\nu)\cdot\delta X\circ X,\theta_0\rangle_2
\end{equation*}
Hence we conclude
\begin{equation*}
 \frac{\delta \sP_1}{\delta\bZ}(X_t)= -\nabla\varphi(X_t^\nu)\theta_0  
\end{equation*}
While for $\sP_2(X)$, denote by $X_\varepsilon^\nu=X_t^\nu+\varepsilon\delta X\circ X_t$ and $A_\varepsilon^\nu=(X_\varepsilon^\nu)^{-1}$, $\mathbf{B}:=\nabla XB_0$, $\tilde B=\mathbf{B}\circ A_t^\nu$ and $B=\mathbb E\tilde B$. Notice $\nabla X=\nabla X^\nu$, we have:
\begin{equation*}
0=\frac{d}{d\varepsilon}\leps (X_\varepsilon^\nu\circ A_\varepsilon^\nu)=\delta X\circ X\circ A^\nu+\nabla X_t^\nu(A_t^\nu)\cdot\partial_\varepsilon A_\varepsilon^\nu    
\end{equation*}
\begin{equation*}
\frac{d}{d\varepsilon}\leps \frac{1}{2}\lVert\mathbb E_W[(\nabla X_\varepsilon B_0)\circ A_\varepsilon^\nu]\rVert_2^2=\mathbb E_W\langle B,  (\nabla(\delta X\circ X)B_0)\circ A_t^\nu\rangle_2+\mathbb E_W\langle B, \nabla\mathbf{B}(A_t^\nu)\cdot (-\nabla A_t^\nu)\delta X(\phi_t)\rangle_2
\end{equation*}
\begin{equation*}
= -\mathbb E_W\langle \nabla\cdot[B(X_t^\nu)\otimes B_0],  \delta X\circ X\rangle_2- \mathbb E_W\langle  (\nabla^*\tilde BB)\circ \phi_t^{-1},\delta X\rangle_2 
\end{equation*}
\begin{equation*}
=-\mathbb E_W\langle B_0\cdot\nabla (B(X_t^\nu))+ (\nabla^*\tilde BB)\circ X_t^\nu,\delta X\circ X\rangle_2     
\end{equation*}
\begin{equation*}
=- \mathbb E_W\langle \nabla X_t^\nu B_0\cdot\nabla B(X_t^\nu)+ (\nabla^*\tilde BB)\circ X_t^\nu,\delta X\circ X\rangle_2  
\end{equation*}
\begin{equation*}
=-\mathbb E_W\langle (\tilde B\cdot\nabla B)\circ X_t^\nu+ (\nabla^*\tilde BB)\circ X_t^\nu,\delta X\circ X\rangle_2.      
\end{equation*}
Hence, we conclude:
\begin{equation*}
 \frac{\delta\sP}{\delta \bZ}(X)=-\nabla\varphi(X_t^\nu)\theta_0-(\tilde B\cdot\nabla B)\circ X_t^\nu-(\nabla^*\tilde BB)\circ X_t^\nu  . 
\end{equation*}
Plugging the above into \eqref{e:Eulerian Momentum}, the critical point momentum equation of the Boussinesq-MHD convection model is:
\begin{equation*}
\partial_t\xi+u\cdot\nabla \xi +\nabla^*u\xi=\nu\Delta\xi+\theta \nabla\varphi+B\cdot\nabla B+\nabla\frac
{1}{2}\lvert B\rvert^2,\quad u=\bP\xi  
\end{equation*}
where $\theta=\mathbb E_W(\theta_0\circ A_t^\nu)$, $B=\mathbb E_W[(X_t^\nu)_\sharp B_0]$ satisfies the following transport-diffusion equation
\begin{equation*}
 \partial_t\theta+u\cdot\nabla\theta=\nu\Delta\theta,\quad \partial_t B+[u,B]=\nu\Delta B .  
\end{equation*}
Thanks to \eqref{e:Dual Lie Transport Projected}, above coupled system is equivalent to viscous Boussinesq equation with MHD convection, in the unit Prandtl number scenario (viscosity, diffusivity and resistivity are the same). The associated stochastic Lagrangian-Eulerian system reads as follows:

\begin{equation}\label{e:Stochastic Lagrangian Boussinesq&MHD}
\left\{
\begin{aligned}
&dX_t^\nu(t,x)=u(t,X_t^\nu(x))dt+\sqrt{2\nu} dW_t,\quad  A_t^\nu(x)=(X_t^\nu)^{-1}(x)\\
&\theta_t=\mathbb E_W[\theta_0\circ A_t^\nu(x)]\\
&B_t=\mathbb E_W \tilde B_t,\quad \tilde B_t=(X_t^\nu)_\sharp B_0\\
&u(t,x)=\mathbf P\mathbb E_W\bigg[\nabla^*A_t^\nu(x)\bigg(u_0+\int_0^t \nabla^*X_\tau^\nu [\theta_0\nabla\varphi(X_\tau^\nu)+ (\tilde B\cdot\nabla B+\nabla^*\tilde BB)(X_\tau^\nu)]d\tau\bigg)\circ A_t^\nu(x)\bigg]
\end{aligned}
\right.
\end{equation}

\paragraph{Example 3.5} Our last example will be the viscous Hall-MHD model. In this example we assume everything to be smooth to avoid issue of rough transport. Recall the following viscous MHD equation with Hall effect ($\alpha>0$):
\begin{equation}\label{e:Hall MHD}
\left\{
\begin{aligned}
&\partial_tu+u\cdot\nabla u+\nabla p=\nu\Delta u+J\times B\\
&\partial_t B+[u,B]-\alpha\nabla\times(\nabla\times B\times B)=\nu\Delta B\\
&\nabla\cdot u=\nabla\cdot B=0\\
&u(0,\cdot)=u_0,\quad B(0,\cdot)=B_0
\end{aligned}
\right.
\end{equation}

Hall MHD can be viewed as a physical system with two natural momentum variables, as mentioned in \cite{Eyi09}. To see this, thanks to $\nabla\cdot u=\nabla\cdot J=0$, apply Lemma 2.4 we pass the magnetic transport equation to the following dual-Lie transport equation of vector potential $\mathbf V:=\BioS B$:
\begin{equation*}
 \partial_t\mathbf V+(u-\alpha J)\cdot\nabla \mathbf V+\nabla^*(u-\alpha J)\mathbf V+\nabla q=\nu\Delta\mathbf{V} :=(I) 
\end{equation*}
where $\nabla q$ is a pressure you gain to imposing constraint $\nabla\cdot \mathbf V=0$. On the other hand we could rewrite the momentum equation of $u$ in \eqref{e:Hall MHD} as
\begin{equation*}
 \partial_tu+u\cdot\nabla u+\nabla^* uu+\underbrace{(\nabla p-\frac{1}{2}\lvert u\rvert^2)}_{:=\nabla\tilde p} =B\cdot\nabla B+\nu\Delta u :=(II).   
\end{equation*}
Hence, computing $(II)-\alpha^{-1}(I)$ one get for $Q=\tilde p-\alpha^{-1}q$ one has
\begin{equation}\label{e:ME1}
 \partial_t(u-\alpha^{-1}\mathbf V)+u\cdot\nabla(u-\alpha^{-1}\mathbf V)+\nabla^*u(u-\alpha^{-1}\mathbf V)+\nabla Q-\nu\Delta (u-\alpha^{-1} \mathbf{V}) =B\cdot\nabla B-(\nabla^* J\mathbf V+ J\nabla \mathbf V).   
\end{equation}

Note that
\begin{equation*}
 \nabla^* J\mathbf V+ J\nabla \mathbf V=\nabla(\mathbf V\cdot J)+ J\times B=B\cdot\nabla B+\nabla(\mathbf V\cdot J-\frac{1}{2}\lvert B\rvert^2).   
\end{equation*}
Therefore, RHS of \eqref{e:ME1} is indeed a gradient and can be absorbed into $Q$. Denote by $\tilde Q$ the new pressure, by setting $v:=u-\alpha J$, we obtain two homogeneous dual Lie transport equations for $\xi:=u-\alpha^{-1}\mathbf V$ and $\eta=\mathbf V$:
\begin{equation*}\label{e:MEHMHD}
\partial_t \xi+u\cdot\nabla\xi+\nabla^*u\xi+\nabla \tilde Q=\nu\Delta\xi,\quad \partial_t\eta+v\cdot\nabla\eta+\nabla^*v\eta+\nabla q=\nu\Delta\eta    
\end{equation*}
which motivates the following two-velocity Lagrangian:
\begin{equation*}
\tilde\Phi(u,v,X,Y,\xi,\eta)
\end{equation*}
\begin{equation*}
=\int_0^T\frac{1}{2}\lVert u\rVert_2^2-\frac{1}{2}\big\lVert\big(\frac{u-v}{\alpha}\big)\big\rVert_{\dot H^{-1}}^2-\mathbb E_W\langle\tilde\xi, \dot X-u(X_t^\nu)\rangle_2-\mathbb E_{\tilde W}\langle\tilde\eta,\dot Y-v(Y_t^\nu)\rangle_2dt    
\end{equation*}
where $W$, $\tilde W$ are two independent Wiener processes. The associated stochastic Lagrangian formulation reads as:
\begin{equation}\label{VHMHDLag}
\left\{
\begin{aligned}
&d X_t^\nu=u(t,X_t^\nu)dt+\sqrt{2\nu}dW_t,\quad A_t^\nu=(X_t^\nu)^{-1}\\
&dY_t^\nu=v(t,Y_t^\nu)+\sqrt{2\nu}d\tilde W_t,\quad Z_t^\nu=(Y_t^\nu)^{-1}\\
&u=\mathbb E_W\mathbf P[\nabla^*A_t^\nu(u_0-\alpha^{-1}{\BioS}B_0)\circ A_t^\nu]+\mathbb E_{\tilde W}\mathbf P[\nabla^*Z_t^\nu(\alpha^{-1}{\BioS}B_0)\circ Z_t^\nu]\\
&v=u-\alpha\nabla\times B\\
&B=\mathbb E_W(Y_t^\nu)_\sharp B_0
\end{aligned}
\right.
\end{equation}

\paragraph{Remark 3.6} Our variational formulation is not satisfactory for deriving general viscous models with diffusive advection, where more than one source of dissipation appears. Such an issue can not be simply overcome by introducing an additional random process to generate the dissipation of the tracer, as discussed in \cite{Ya18}. On the other hand, coupled fluid-tracer models with arbitrary viscosity and  diffusivity can be produced under metricplectic framework, embracing a more geometric interpretation of diffusion in the sense of Wasserstein geometry \cite{JKO98}. Extensive discussion of this aspect can be found in a companion paper \cite{MP26}.

\subsection{A Generalized Flow Viewpoint}

In \cite{Iy06} the author also discussed a non-averaged modification of \eqref{e:Lagrangian NSE}. In that scenario, due to lack of averaging and smoothing mechanism, the system does not produce solution of Navier-Stokes, but gives a random perturbation of Euler flow instead. It is natural to consider such models in the above variational framework and we show that they can be derived via a variational principle which reminiscent of Brenier's generalized (or multi-phased) flow \cite{Br89}.

We first briefly recall the notion of generalized incompressible flow (See \cite{Sh00}). Fix a probability space $(\Omega,\mathcal A,\mathbb P)$ and a random probability measure $\rho_0^\omega(dx)$ on $D$, a generalized incompressible flow is a random measurable map 
\begin{equation*}
X:(t,\omega,x)\in [0,T]\times\Omega\times D\to X_t^\omega(x)\in D    
\end{equation*}
such that for any smooth functions $f\in C_c^\infty(D)$:
\begin{equation*}
\int_{\Omega} \int_D f(X_t^\omega(x))\rho_0^\omega(dx)\mathbb P(d\omega)=\int_D f(x)dx  
\end{equation*}
that is, the following incompressibility condition holds: Denote by $\mathcal L_D$ the normalized Lebesgue measure confined to $D$, we have:
\begin{equation*}
\int_\Omega (X_t^\omega)_\sharp\rho_0^\omega(dx)\mathbb P(d\omega)=\mathcal L_D(dx)
\end{equation*}
Notice here we use the product space $\Omega\times D$ in place of generic probability space $\Omega$ in \cite{Sh00} to keep consistency with our variational framework. 

We now formally discuss the connection between the stochastic Hamilton-Pontryagin principle and generalized flow. The idea is to send $\nu=0$ but still keep the randomness survive. To make sense of such procedure, we randomize our velocity field also and we relax the incompressibility to an averaged sense, instead of pathwise (like $u_\nu$).

Introduce the following Lagrangian of random elements:
\begin{equation*}
\sL(u^\omega,X^\omega,\xi^\omega)=\int_0^T\int_\Omega\int_D \frac{1}{2}\lvert u_t^\omega(X_t^\omega(x))\rvert^2+\xi_t^\omega(x)\cdot [\dot X^\omega-u_t^\omega(X_t^\omega(x))]\rho_0^\omega(dx)\mathbb P(d\omega)dt
\end{equation*}

To fix idea, here we assume sufficient spatial regularity. Now the following theorem holds:
 \paragraph{Theorem 3.7} Under the incompressibility constraint:
\begin{equation}
\int_{\Omega} (X_t^\omega)_\sharp\rho_0(dx)\mathbb P(d\omega)=\mathcal L_D(dx)
\end{equation} 
the following statements hold:
\begin{itemize}
    \item (i) Assume $(u,X,\xi)$ is a critical point equation of $\sL$, then $(u^\omega,\rho^\omega):=(u,X_\sharp\rho_0)$ solves the following generalized Euler's equation
 \begin{equation}\label{Generalized Euler}
\left\{
\begin{aligned}
&\partial_t (\rho_t^\omega u_t^\omega)+\nabla\cdot(\rho_t^\omega u_t^\omega\otimes u_t^\omega)+\rho_t^\omega\nabla p=0\\
&\partial_t\rho_t^\omega+\nabla_x\cdot(\rho_t^\omega u_t^\omega)=0
\end{aligned}
\right.
\end{equation}

\item (ii) Consider Hamiltonian:
\begin{equation*}
\tilde\sH(\rho_t^\omega,\chi_t^\omega,u_t^\omega)=\int_0^T\int_\Omega\int_{\mathbb R^d} -\frac{1}{2}\lvert u_t^\omega(x)\rvert^2\rho_t^\omega(x)-\chi_t^\omega\nabla\cdot(\rho_t^\omega u_t^\omega) dx\mathbb P(d\omega)dt
\end{equation*}

Assume $(\rho^\omega,u^\omega,\chi^\omega)$ is a triplet of solution of the following Hamiltonian ODE
\begin{equation}\label{e:HamODE2}
\left\{
\begin{aligned}
&\partial_t\rho_t^\omega=\frac{\delta \tilde\sH}{\delta \chi}\\
&{\partial_t \chi_t^\omega}=-\frac{\delta\tilde\sH}{\delta\rho}\\
&\frac{\delta \tilde\sH}{\delta u}=0
\end{aligned}
\right.
\end{equation}
then $(u^\omega,\rho^\omega)$ solves the generalized Euler equation.

\end{itemize}
 
 \textbf{Proof.} (i) We compute the first variation of $\sL$. The variation  $\delta u$ and $\delta\xi$ are free, thus arbitrary elements in $[0,T]\times\Omega\times D\to \mathbb R^d$. While for $\delta X$ we consider compactly supported variations $\delta X$. Denote by $X_\varepsilon=X+\varepsilon\delta X$, we therefore have for any test function $f\in C_c^\infty(D)$:
 \begin{equation}
 \frac{d}{d\varepsilon}\leps\int_{\Omega\times D} f(X_\varepsilon)\rho_0^\omega (dx)\mathbb P(d\omega)=\int_\Omega \int_{D} \nabla f(X_t^\omega(x))\delta X_t^\omega(x)\rho_0^\omega(dx)\mathbb P(d\omega)=0
 \end{equation}
 Now compute variational derivative of each variables $u$, $\xi$, $X$ and taking above flexibility condition for $\delta X$ into account, we have
 \begin{equation*}
\left\{
\begin{aligned}
&\dot X_t^\omega(x)=u_t^\omega(X_t^\omega(x))\\
&\dot \xi_t^\omega(x)+\nabla^*u_t^\omega(X_t^\omega(x))\xi_t^\omega(x)-\frac{1}{2}\nabla\lvert u_t^\omega\rvert^2 (X_t^\omega(x))=\nabla p(X_t^\omega(x))\\
& \xi_t^\omega=u_t^\omega\circ X_t^\omega
\end{aligned}
\right.
\end{equation*}

Therefore plug the first and third equation into the second we conclude
\begin{equation}\label{Multiphase Momentum}
(\partial_tu_t^\omega+u_t^\omega\cdot\nabla u_t^\omega)(X_t^\omega(x))=\nabla p(X_t^\omega(x)).
\end{equation}
Use the notation $\rho_t^\omega=(X_t^\omega)_\sharp \rho_0^\omega$, we see that $\rho_t^\omega$ solves the following continuity equation:
\begin{equation}\label{Multiphase Density}
\partial_t\rho_t^\omega+\nabla_x\cdot(\rho_t^\omega u_t^\omega)=0.
\end{equation}
Now thanks to \eqref{Multiphase Momentum}, \eqref{Multiphase Density} it's straightforward to verify:
\begin{equation*}
\partial_t(\rho_t^\omega u_t^\omega)=-\nabla\cdot(\rho_t^\omega u_t^\omega)u_t^\omega+\rho_t^\omega(\nabla p-u_t^\omega\cdot\nabla u_t^\omega)=-\nabla\cdot(\rho_t^\omega u_t^\omega\otimes u_t^\omega)+\rho_t^\omega\nabla p.
\end{equation*}
Which is the generalized Euler equation.

\quad

(ii) Compute variation with respect to each variable, here variation in $\chi$ and $u$ are free, while variation in $\rho$ satisfies
\begin{equation*}
\int_{\Omega}\int_Df(x)\delta\rho_t^\omega(dx)\mathbb P(d\omega)=0,\quad\text{ for all }f\in C_c(D)    .
\end{equation*}
Therefore the variational equation reads as:
\begin{equation*}
\frac{\delta \tilde\sH}{\delta \chi}=-\nabla\cdot(\rho u),\quad \frac{\delta\tilde \sH}{\delta \rho}=u\cdot\nabla\chi-\frac{1}{2}\lvert u\rvert^2+p,\quad \frac{\delta \tilde\sH}{\delta u}=-u\rho+\rho\nabla\chi    
\end{equation*}
where $p$ is independent of $\omega$.   
Now the first equation is simply the continuity equation, while the third equation reads $u_t^\omega=\nabla\chi_t^\omega$ guarantees that $u$ is pathwisely gradient. Now putting above identity into the second equation we have:
\begin{equation*}
\partial_t\chi_t^\omega+\frac{1}{2}\lvert\nabla\chi_t^\omega\rvert^2+ p=0 .
\end{equation*}
To simplify notation we drop the $\omega$ superscript for randomness, now straightforward computation yields:
\begin{equation*}
\partial_t(\rho u)=-u\nabla\cdot(\rho u) -(\nabla p+\nabla^* uu)\rho=-\nabla\cdot(\rho u\otimes u)-\rho\nabla p 
\end{equation*}
where $\rho$, $u$ depends on $\omega$ while $p$ is $\omega$-independent. Therefore we recover \eqref{Generalized Euler} and the proof is completed. \qed

\quad

Above Hamiltonian formulation is a randomized version of Clebsch principle for fluids (cf. \cite{HMR98}\cite{MW83}) and also mirrors the groupoid Hamiltonian formulation of generalized flows in \cite{IK24}.

We now turn to study of the non-averaged model with Brownian perturbation. To this end we fix our probability space $(\Omega,\mathbb P)$ to be Wiener space $(\Gamma_D,\mathcal W)$, which is the canonical space for Brownian motion. Here $\mathcal W$ is the Wiener measure on $\Gamma_D:=C([0,T];D)$. Now denote by
\begin{equation*}
 X_{t,\nu}^\omega=X_t^\omega+\sqrt{2\nu}\omega(t).   
\end{equation*}
Consider the perturbed Hamiltonian:
\begin{equation}
\tilde\sH_\nu(u^\omega,X^\omega,\xi^\omega)=\int_0^T\int_\Omega\int_D\big[-\frac{1}{2}\lvert u_t^\omega(X_{t,\nu}^\omega(x))\rvert^2+\xi_t^\omega(x)u_t^\omega(X_{t,\nu}^\omega(x))\big]\rho_{0}^{\omega}(dx)\mathbb P(d\omega)dt.
\end{equation}
We have the following 
\paragraph{Theorem 3.8}  Assume $(u,X,\xi)$ solves the Hamiltonian equation associated with $\tilde\sH_\nu$, that is:
\begin{equation}\label{e:HamODE3}
\left\{
\begin{aligned}
&\dot\xi_t^\omega=-\frac{\delta \tilde\sH_\nu}{\delta X}\\
&\dot X_t^\omega=-\frac{\delta\tilde\sH_\nu}{\delta\xi}\\
&\frac{\delta \tilde\sH_\nu}{\delta u}=0
\end{aligned}
\right.
\end{equation}

Then, under initial condition $X_0(x)=x$, $u^\omega(0,x)=\bar u$, one has $u^\omega(t,x)=u_E(x-\sqrt{2\nu} W_t^\omega)$ for a.s. $\omega$, where $u_E$ is the Euler flow associated with initial data $u(0,\cdot)=\bar u$.

\textbf{Proof.} Similar computation as in theorem 3.7 yields:
\begin{equation*}
 \frac{\delta \tilde\sH_\nu}{\delta X}=\nabla^*u(X_{t,\nu}^\omega)\xi-\frac{1}{2}\nabla\lvert u\rvert^2(X_{t,\nu}^\omega)+\nabla p(X_t^\omega).   
\end{equation*}
Therefore the Hamiltonian ODE \eqref{e:HamODE3} reads as
\begin{equation*}
\left\{
\begin{aligned}
&\dot\xi_t^\omega+\nabla^*u(X_{t,\nu}^\omega)\xi_t^\omega-(\nabla^* uu)(X_{t,\nu}^\omega)+\nabla p(X_t^\omega)=0\\
&\dot X_t^\omega=u(X_{t,\nu}^\omega)\\
&\xi_t^\omega=u(X_{t,\nu}^\omega)
\end{aligned}
\right.
\end{equation*}
Notice the second equation implies that 
\begin{equation*}
 dX_{t,\nu}^\omega=u_t^\omega(X_{t,\nu}^\omega)dt+\sqrt{2\nu}dW_t^\omega.   
\end{equation*}
Therefore by Ito-Wenztell formula one conclude $u_t^\omega=\xi_t^\omega\circ A_{t,\nu}^\omega$ satisfies:
\begin{equation*}
du_t^\omega+(u_t^\omega\cdot\nabla u_t^\omega-\nu\Delta u_t^\omega+\nabla p_\nu(x))dt+\sqrt{2\nu}\nabla u_t^\omega dW_t^\omega=0,\quad u_0^\omega(x)=\bar u.    
\end{equation*}
where $p_\nu(x)=p(x-\sqrt{2\nu}W_t)$.

Now assume $u_E(t)$ is the Euler flow starting from initial data $\bar u$, by Ito's formula the random shift $u_E^\omega(t,x):=u_E(x-\sqrt{2\nu} W_t^\omega)$ satisfies the same equation as $u_t^\omega$. Hence uniqueness of solution of above SPDE guarantees $u_t^\omega(x)=u_E(t,x-\sqrt{2\nu}W_t^\omega)$. \qed

\subsection{A Generalized Circulation Theorem}

We close this section with discussions on conserved quantities due to particle relablling symmetry (by virtue of Noether's theorem) from our Lagrangian formulation, in both inviscid and viscous scenarios. Indeed all such conservation laws are known as Casimir invariants in Language of geometric mechanics, referring to constants of coadjoint motions, which could be regarded as a simple duality result in our formulation. Its generalization to viscous fluids was observed in \cite{Dr22} for Navier-Stokes. Here we first state an inviscid version with generalization to many other fluid models.

Let $v_t$ be a weak solution of \eqref{e:Lie Transport} with initial data $v_0$, we typically highlight the case $v\in \mathcal M_\sigma$, i.e. $v$ is a div-free vector-valued measure. Here we let $\mathcal M(D;\mathbb R^d)$ be the space of $\mathbb R^d$-valued measures on $D$ and define:
\begin{equation*}
\mathcal M_\sigma(D)=\{\mu\in \mathcal M(D;\mathbb R^d): \int_D \varphi\cdot d\mu=0\text{ for all }\varphi\in C_c^\infty(D;\mathbb R^d)\} .   
\end{equation*}
In the case $d=2$, one could think of $\mu=\nabla^\perp f$ in distributional sense for $f\in L^1(D)$. Now, given flow $X_t$ of $u\in H_\sigma^s$, we define the push-forward of $v$ by $X_t$ using the following duality identity:
\begin{equation}\label{Pushforward of vector Measure}
\langle \psi,(X_t)_\sharp v\rangle_*:=\langle (X_t)^*\psi,v\rangle_* \quad\text{for all }\psi\in C_c^\infty(D;\mathbb R^d).   
\end{equation}
Here 
\begin{equation*}
(X_t)^* \psi=\nabla^*X_t \psi\circ X_t    
\end{equation*}
is the pull-back of vector field $\psi$ in the sense of co-vector. By above duality definition,  distributional solution of \eqref{e:Lie Transport} with div-free measure initial data $v$ is given by $v_t:=(X_t)_\sharp v$. Moreover $v_t\in \mathcal M_\sigma$ since for any $\varphi\in C_c^\infty(D)$:
\begin{equation*}
\langle (X_t)_\sharp v,\nabla\varphi\rangle_*=\langle v,\nabla^*X_t\nabla\varphi\circ X_t\rangle_*=\langle v,\nabla(\varphi\circ X_t)\rangle_*=0.    
\end{equation*}
We now state:

\paragraph{Theorem 3.9(Generalized Kelvin Theorem, cf.\cite{Dr22})} Assume $w_t$ is a strong solution of dual Lie transport equation \eqref{e:Dual Lie Transport} with drift $u\in C([0,T];H_\sigma^s)$ and initial data $w_0\in H^{s-1}$, then for $v_t$ which is a weak solution of \eqref{e:Lie Transport} with initial data $v\in \mathcal M_\sigma$, we have:
\begin{equation}\label{e:Generalized Kelvin}
 \langle w_t,v_t\rangle_*=\langle w_0,v\rangle_* .  
\end{equation}

Proof: We provide a one-line Lagrangian proof here. By Lemma 2.4 we have:
\begin{align*}
\langle w_t,v_t\rangle_*=\langle (A_t)^*w_0,(X_t)_\sharp v_0\rangle_*=\langle (X_t)^*(A_t^*w_0),v_0\rangle_*=\langle w_0,v_0\rangle_*. \qed
\end{align*}
\paragraph{Remark 3.10} One may directly obtain a generalization of the  following Ertel's theorem: Let $w_t$ be strong solution of \eqref{e:Dual Lie Transport} with drift $u\in C([0,T];H_\sigma^s)$ and let $v_t$ be a strong solution of \eqref{e:Lie Transport}, with initial data $w_0$, $v_0$ respectively, then $\lambda_t=w_t\cdot v_t$ is a constant of motion. To see this, notice in Lagrangian variables one has:
\begin{equation*}
w_t\cdot v_t=\nabla^*A_tw_0\circ A_t\cdot (\nabla X_tv_0)\circ A_t=\nabla^*X_t\circ A_t\nabla^*A_tw_0\circ A_t\cdot v_0\circ A_t=(w_0\cdot v_0)\circ A_t    
\end{equation*}
 Denote by $\Gamma$ the space of all simple closed curves on $\mathbb R^d$ with smooth parametrization on $[0,1]$, notice that we could identify circulation along loop as a vector-valued measure ($1$-current): for any smooth divergence free vector field $\psi$ we have
\begin{equation*}
 \oint_\gamma \psi\cdot dl=\langle \psi,\mu_\gamma\rangle_*    
\end{equation*}
where above paring should be considered as duality pairing between vector field and vector-valued measure(See \cite{ABC13}). Measure $\mu_\gamma$ in explicit form is given by
\begin{equation}
 \mu_\gamma=\vec{\tau_\gamma}\mathscr H_{C_\gamma}^{1}(dx),\quad C_\gamma=\mathsf{Im}(\gamma)\subset\mathbb R^d   
\end{equation}
here $\vec{\tau_\gamma}$ is the unit tangent vector of curve $\gamma$. 

To convince the reader that theorem 3.9 is a generalized Kelvin's theorem, it suffices to show that 
 for measure-preserving flow $X$ of $u$ the following holds:
\begin{equation*}
\mu_{X_t(\gamma)}=(X_t)_\sharp(\mu_\gamma)    
\end{equation*}
To see this, choose any smooth test vector field $\psi$ and we observe
\begin{equation*}
 \langle\mu_{X_t(\gamma)},\psi\rangle_*=\oint_{X_t(\gamma)}\psi\cdot dl = \int_0^1 \psi(X_t\circ\gamma)(s)\cdot\frac{d}{ds}(X_t(\gamma(s)))ds 
\end{equation*}
\begin{equation*}
 =\int_0^1 (\nabla^*X_t\psi\circ X_t)(\gamma(s))\cdot\dot\gamma(s)ds   
\end{equation*}
\begin{equation*}
 =\oint_\gamma (\nabla^*X\psi\circ X)\cdot dl=\langle \mu_\gamma,(X_t)^*\psi\rangle_*=\langle (X_t)_\sharp\mu_\gamma,\psi\rangle_*  
\end{equation*}
We remark that theorem 3.9 recovers conservation of circulation/helicity in various inviscid incompressible hydrodynamic models.

Now we move on to the viscous setting and introduce stochastic counterpart of generalized circulation theorem \eqref{e:Generalized Kelvin}, which generalizes the pathwise Kelvin theorem and statistical Kelvin theorem discussed in \cite{CI08}\cite{DH18}\cite{Eyi09}.

\paragraph{Theorem 3.11} Let $\sP=0$, assume $(u,X,\tilde\xi)$ is a critical point of \eqref{e:SFunc}, denote by $\xi^\prime:=\tilde\xi\circ A_t^\nu$, then the following holds:
\begin{itemize}
    \item (I) (Pathwise Kelvin) Consider the random div-free vector measure $\tilde v_t=(X_t^\nu)_\sharp v$ with $v\in \mathcal M_\sigma$, then we have the following pathwise conservation:
    \begin{equation}\label{Pathwise Kelvin}
     \langle \xi_t^\prime,\tilde v_t\rangle_*(\omega)=\langle \xi_0,v\rangle_*\quad\text{for a.e. }\omega\in \Omega, t\in [0,T]   
    \end{equation}
    \item (II) (Statistical Kelvin) For any $\mu\in \mathcal M_\sigma$, we have
    \begin{equation}\label{Statistical Kelvin}
     \langle \xi_t,\mu\rangle_*=\langle \xi_0,\mathbb E_W(A_t^\nu)_\sharp \mu\rangle_*   
    \end{equation} 
\end{itemize}

Proof: The proof of (I) follows from straightforward computation. Notice that
\begin{equation*}
\xi_t^\prime=\tilde\xi_t\circ A_t^\nu=[(\nabla^* X_t^\nu)^{-1}\xi_0]\circ A_t^\nu=\nabla^*A_t^\nu\xi_0\circ A_t^\nu .   
\end{equation*}
Therefore we have:
\begin{equation*}
 \langle \xi^\prime_t,\tilde v_t\rangle_*=\langle(X_t^\nu)^*\xi^\prime_t,  v\rangle_*=\langle(X_t^\nu)^*((A_t^\nu)^*\xi_0),v\rangle_*=\langle (X_t^\nu\circ A_t^\nu)^* \xi_0,v\rangle_*=\langle \xi_0,v\rangle_* \text{ for a.e. }\omega .  
\end{equation*}

(III) simply follows from theorem 2.6 and Fubini:
\begin{equation*}
\langle \xi_t, \mu\rangle_*=\mathbb E_W\langle (A_t^\nu)^*\xi_0,\mu\rangle_*=\langle\xi_0,\mathbb E(A_t^\nu)_\sharp \mu\rangle_*. \qed    
\end{equation*}

We now illustrate how above theorem yields various conservation laws stemming from particle relabeling symmetry in hydrodynamic setting.

\paragraph{Example 3.12} (I) Assume $d=3$, $\bomega^\prime:=\nabla\times \xi^\prime$, hence $\bomega=\mathbb E_W\bomega^\prime$. Then $\bomega^\prime=(X_t^\nu)_\sharp \bomega$. To see this, recall that for $\dot X=u_\nu(X)$, we have $A_t^\nu= A_t\circ \phi_t$, therefore by Lemma 2.5, pathwisely we have:
\begin{equation*}
\bomega^\prime=\nabla\times [(A_t\circ \phi_t)^*\xi_0]=\nabla\times [\phi_t^*(A_t)^*\xi_0]=[\nabla\times (A_t^*\xi_0)]\circ \phi_t  
\end{equation*}
\begin{equation*}
=  [(X_t)_\sharp\bomega_0]\circ \phi_t=(\phi_t^{-1})_\sharp [(X_t)_\sharp \bomega_0]=(\phi_t^{-1}\circ X_t)_\sharp \bomega_0=(X_t^\nu)_\sharp \bomega_0  
\end{equation*}
Hence by Theorem 3.11 (I) we have the following pathwise conservation of helicity:
\begin{equation*}
 \langle \xi_t^\prime,\bomega_t^\prime\rangle_2= \langle \xi_0,\bomega_0\rangle_2=\langle \underbrace{\bP\xi_0}_{=u_0},\bomega_0\rangle  
\end{equation*}
In case $d=2$, we have $\bomega^\prime=\bomega_0\circ A_t^\nu$ is a scalar, and $\nabla^\perp \bomega^\prime$ satisfies:
\begin{equation*}
\nabla^\perp\bomega^\prime=\bJ \nabla(\bomega_0\circ A_t^\nu)=\bJ \nabla^*A_t^\nu\nabla\bomega_0\circ A_t^\nu=(\nabla A_t^\nu)^{-1}\bJ \nabla\bomega_0(A_t^\nu)=(X_t^\nu)_\sharp\nabla^\perp\bomega_0.     
\end{equation*}
Hence above pathwise conservation law reduces to the vorticity conservation
\begin{equation*}
\lVert\bomega_0\rVert_2^2 =\langle \xi_0,\nabla^\perp\bomega_0\rangle_2= \langle \xi_t^\prime,\nabla^\perp\bomega_t^\prime\rangle_2=\langle \nabla\times \xi_t^\prime,\bomega_t^\prime\rangle_2=\lVert\bomega_t^\prime\rVert_2^2 . 
\end{equation*}

On the other hand for MHD, let $d=3$ we have pathwise conservation of magnetic helicity: Let $\tilde A=\BioS \tilde B=(A_t^\nu)^*(\BioS B_0)$, then we have
\begin{equation*}
\langle \tilde A_t,\tilde B_t\rangle_2=\langle (A_t^\nu)^*A_0,(X_t^\nu)_\sharp B_0\rangle_2=\langle A_0,B_0\rangle_2 .   
\end{equation*}
While in general dimension the following conservation of cross helicity holds:
\begin{equation*}
 \langle \tilde \xi_t,\tilde B_t\rangle_2=\langle (A_t^\nu)^*\xi_0,(X_t^\nu)_*B_0\rangle_2=\langle \xi_0,B_0\rangle_2  . 
\end{equation*}
(II) Choose $v=\mu_\gamma$ as in \eqref{Pathwise Kelvin}, we obtain the pathwise conservation of circulation in various viscous hydrodynamic models. For instance, in the case of generalized SQG (Example 3.3), we have for $\tilde v_t=\mu_{X_t^\nu(\gamma)}$ and $\tilde u_t=\bP\xi_t^\prime$:
\begin{equation*}
\oint_{X_t^\nu(\gamma)} \xi_t^\prime\cdot dl=\oint_{X_t^\nu(\gamma)} (-\Delta)^{-\alpha}\tilde u_t\cdot dl=\oint_\gamma (-\Delta)^{-\alpha}u_0\cdot dl   .
\end{equation*}
While in the case of Hall MHD (Example 3.5), let $\xi^\prime:=\tilde\xi\circ A_t^\nu$ and $\eta^\prime:=\tilde\eta\circ Z_t^\nu$. By Lemma 2.4, we have:
\begin{equation*}
 \eta^\prime=\tilde\eta\circ Z_t^\nu=\nabla^*Z_t^\nu(\BioS B_0)\circ Z_t^\nu=\BioS[(Y_t^\nu)_\sharp B_0] .  
\end{equation*}
Hence, we have
\begin{equation*}
\oint_{X_t^\nu(\gamma)} \xi_t^\prime\cdot dl=\oint_\gamma (\tilde u_t-\alpha^{-1}\BioS\tilde B_t)\cdot dl=\oint_\gamma (u_0-\alpha^{-1}\BioS B_0)\cdot dl .  
\end{equation*}
and
\begin{equation*}
\oint_{Y_t^\nu(\gamma)}\eta^\prime \cdot dl=\oint_{Y_t^\nu(\gamma)} \BioS \tilde B_t\cdot dl=\oint_\gamma\BioS B_0 \cdot dl=\oint_\gamma \eta_0 \cdot dl.   
\end{equation*}

Notice that applying Stokes theorem to the above circulation conservation we obtain the following stochastic Alfv\'en theorem: For any smooth surface $S\subset \mathbb R^d$ we have:
\begin{equation*}
 \int_{Y_t^\nu(S)} \tilde B_t\cdot d\mathbf{S}=\int_{S}B_0\cdot d\mathbf{S}.   
\end{equation*}

(III) Assume moreover $v\in L_\sigma^2(D;\mathbb R^d)$, pointwisely we have:
\begin{equation*}
\xi_t^\prime\cdot \tilde v_t=(\nabla^* A_t^\nu\xi_0\circ A_t^\nu,(\nabla X_t^\nu v)\circ A_t^\nu)=(\nabla^* X_t^\nu\circ A_t^\nu\nabla^*A_t^\nu \xi_0\circ A_t^\nu, v\circ A_t^\nu)=(\xi_0,v)\circ A_t^\nu   
\end{equation*}
Therefore let $\lambda_t:=\mathbb E_W(\xi_t^\prime\cdot \tilde v_t)$, we have $\lambda_t=\mathbb E_W[(\xi_0,v)\circ A_t^\nu]$, hence we prove a generalization of Ertel's theorem to viscous setting since $\lambda_t$ satisfies the following parabolic PDE:
    \begin{equation*}
      \partial_t \lambda+u\cdot\nabla\lambda=\nu\Delta\lambda,\quad\lambda(0,\cdot)=\xi_0\cdot v
    \end{equation*}

(IV) Now fix any smooth loop $\gamma\in \Gamma$, smooth surface $S$ and $\mu=\mu_\gamma\in \mathcal M_\sigma$, in case of generalized SQG we have the following statistical Kelvin's theorem.
\begin{equation*}
 \mathbb E_W\oint_{A_t^\nu(\gamma)}(-\Delta)^{-\alpha} u_0\cdot dl=\oint_{\gamma}(-\Delta)^{-\alpha}u_t\cdot dl   
\end{equation*}
  While for Hall-MHD the following holds:
 \begin{equation*}
\mathbb E_W\oint_{A_t^\nu(\gamma)} (u_0-\alpha^{-1}\BioS B_0)\cdot dl=\oint_\gamma v_t\cdot dl,\quad \mathbb E_W\int_{Z_t^\nu(S)} B_0\cdot d\mathbf{S}=\int_S B_t \cdot d\mathbf{S}.       
 \end{equation*} 

 \paragraph{Remark 3.13} In \cite{Dr22}, the author derive a generalized Kelvin's theorem from Eulerian viewpoint, as a generalization of theorem 3.10 to viscous case by adding forward/backward diffusion to the pair of dual transport equations. In Lagrangian viewpoint this resonants the fact that statistics of tangent vectors supported on Lagrangian loops along backward stochastic flows solves the backward Lie transport-diffusion equation.

\section{Local Well-Posedness via Lagrangian Formulation}

We now attempt to obtain local well-posedness result for PDEs via the associated Lagrangian system derived from variational principle. The basic strategy is to use Picard iteration to find the fixed point. We'll establish a self-contained local existence result in Sobolev spaces. The idea is similar to \cite{Iy06} \cite{Zh10}, where the Navier-Stokes case is treated in $C^{k,\alpha}$ and $W^{k,p}$ setting. 

We'll show local well-posedness in $H^s$ of a Boussinesq-MHD system (our Example 3.5), which can be derived from (deterministic) Hamilton-Pontryagin principle with Lagrangian:
\begin{equation*}
\sE(v)=\frac{1}{2}\lVert v\rVert_2^2,\quad \sP(\bZ)=-\langle \varphi,\bZ_\sharp\theta_0\rangle_2+\frac{1}{2}\lVert \nabla\bZ B_0\rVert_2^2.   \end{equation*}
Here we forget about randomness (let $\nu=0$) and well-posedness for corresponding viscous fluids could be shown $mutanis$-$mutandis$, as proof in \cite{Iy06} suggests.

\paragraph{Theorem 4.1} The Lagrangian system \eqref{e:Stochastic Lagrangian Boussinesq&MHD} with $\nu=0$ and $\varphi(x)=x_d$ admits a unique fixed-point solution $(u,X)$ up to time $T$ for $T<C(s,u_0,\theta_0,B_0)$ with 
\begin{equation*}
 C(s,u_0,\theta_0,B_0)\sim C(s)(1+\lVert u_0\rVert_{H^s}^2+\lVert \theta_0\rVert_{H^s}^2+\lVert B_0\rVert_{H^s}^2)^{-1/2}  
\end{equation*}
\textbf{Proof.} We divide our proof into several steps.
\begin{itemize}
    \item \textbf{Step 1.} We first pass to vorticity formulation to obtain desired estimates. Thanks to lemma 2.4 we conclude $\omega=\nabla\times u$ satisfies:
    \begin{equation}\label{e:Vorticity equation}
      \omega(t,x)=(X_t)_\sharp \bigg[\omega_0+\int_0^t (A_\tau)_\sharp [\nabla\times\tilde f(\theta,B)]d\tau\bigg)\bigg]  
    \end{equation}
where 
\begin{equation*}
 \nabla\times\tilde f(\theta,B)=\nabla\times\big(\nabla\theta x_d\big)+\mathscr L_BJ=\nabla\theta\times e_d+\mathscr L_BJ,\quad J=\nabla\times B   
\end{equation*}
Also consider the vorticity form of Lie-transport equation of $B$: By lemma 2.5 we have
\begin{equation}\label{e:Vorticity MHD}
J_t=(X_t)_\sharp\bigg[J_0+\int_0^t (A_\tau)_\sharp(\mathscr L_{B_\tau}\omega_\tau)d\tau+\int_0^t (A_\tau)_\sharp \nabla\times M_\tau d\tau\bigg].    
\end{equation}
Here $J_0=\nabla\times B_0$ and $\nabla\times M$ is a combination of $\partial_iu\partial_jB$ which admits the following bound thanks to $s-1>d/2$ and algebra structure of $H^{s-1}$:
\begin{equation*}
\lVert \nabla\times M_\tau\rVert_{H^{s-1}}\le C\lVert \nabla u_\tau\rVert_{H^{s-1}}\lVert \nabla B_\tau\rVert_{H^{s-1}}\le C\lVert \omega_\tau\rVert_{H^{s-1}}\lVert J_\tau\rVert_{H^{s-1}}.    
\end{equation*}
Therefore sum up \eqref{e:Vorticity MHD} and \eqref{e:Vorticity equation} we obtain the following Lagrangian representation for $\omega+J$:
 \begin{equation*}
\omega_t+J_t=(X_t)_\sharp\bigg[(\omega_0+J_0)+\int_0^t (A_\tau)_\sharp (\mathscr L_{B_\tau}(\omega_\tau+J_\tau))+\int_0^t(A_\tau)_\sharp \nabla\times\big(\nabla\theta x_d+M_\tau\big)d\tau\bigg].     
 \end{equation*}  
 while for $\omega_t-J_t$ we similarly have:
\begin{equation*}
\omega_t-J_t=(X_t)_\sharp\bigg[(\omega_0-J_0)+\int_0^t (A_\tau)_\sharp (\mathscr L_{B_\tau}(\omega_\tau-J_\tau))+\int_0^t(A_\tau)_\sharp \nabla\times\big(\nabla\theta x_d-M_\tau\big)d\tau\bigg].     
\end{equation*} 
Now since push-forward commutes with Lie-bracket we have:
\begin{equation*}
(A_\tau)_\sharp (\mathscr L_{B_\tau}(\omega_\tau-J_\tau))=\mathscr L_{(A_\tau)_\sharp {B_\tau}}[(A_\tau)_\sharp(\omega_\tau-J_\tau)]=\mathscr L_{B_0}[(A_\tau)_\sharp(\omega_\tau-J_\tau)].    
\end{equation*}
Since $B_t=(X_t)_\sharp B_0$.

Denote by $(A_t)_\sharp(\omega_t+J_t)=\zeta_t^+$, $(A_t)_\sharp(\omega_t-J_t)=\zeta_t^-$. Then the following estimate holds:
\begin{equation*}
\lVert\omega_t\pm J_t\rVert_{H^{s-1}}\le \lVert(\nabla X_t\zeta_t^{\pm})\circ A_t\rVert_{H^{s-1}}\le [1+C\exp(CTU)]^{\lfloor s\rfloor+1}\lVert \zeta_t^{\pm}\rVert_{H^{s-1}} .    
\end{equation*}
Moreover, the following polarization inequality holds:
\begin{equation*}
\frac{1}{2}\lVert\omega\rVert_{H^{s-1}}^2+\frac{1}{2}\lVert J\rVert_{H^{s-1}}^2\le \frac{1}{4}\big(\lVert\omega_t+ J_t\rVert_{H^{s-1}}^2+\lVert\omega_t- J_t\rVert_{H^{s-1}}^2\big)   
\end{equation*}
\begin{equation} \label{e:EnergyEs}
\le \frac{1}{4}[1+C\exp(CTU)]^{\lfloor s\rfloor+1}(\lVert \zeta_t^{+}\rVert_{H^{s-1}}^2+\lVert \zeta_t^{-}\rVert_{H^{s-1}}^2).    
\end{equation}
\item \textbf{Step 2.} Now we establish estimate for $\zeta^{\pm}$. Observe that $\zeta_t^{\pm}$ satisfies following equation:
\begin{equation*}
\zeta_t^{\pm}= \zeta_0^{\pm}+\int_0^t \mathscr L_{B_0}\zeta_\tau^{\pm} d\tau+\int_0^t (A_\tau)_\sharp\big(\nabla\theta_\tau\times e_d \pm \nabla \times M_\tau\big)d\tau .
\end{equation*}
Thus $\zeta_t^{\pm}$ satisfies the  Lie transport equation with vector field $-B_0\in H_\sigma^s$ and  force
 \begin{equation*}
 \nabla\times f_t^{\pm}=(A_t)_\sharp\big(\nabla\theta_t\times e_d \pm \nabla \times M_t\big)    
 \end{equation*}
 Assume $\lVert \nabla B_0\rVert_{H^{s-1}}\le U$ , which is the case we'll work with later. By \eqref{e:Lie transport Estimate} we have
 \begin{equation*}
 \lVert \nabla\times f^{\pm}\rVert_{C_tH_x^{s-1}}\le C[1+C\exp(CTU)]^{\lfloor s\rfloor+1}\lVert R^{\pm}\rVert_{C_tH_x^{s-1}} 
 \end{equation*}
 where term $R$ satisfies the following estimate:
 \begin{equation*}
 \lVert R^{\pm}\rVert_{C_tH_x^{s-1}}=\big\lVert\nabla\theta_t\times e_d \pm \nabla \times M_t\big\rVert_{C_tH_x^{s-1}}    
 \end{equation*}
 \begin{equation*}
 \le C\lVert\nabla\theta\rVert_{C_tH_x^{s-1}}+C\lVert\omega\rVert_{C_tH_x^{s-1}}\lVert J\rVert_{C_tH_x^{s-1}}     
 \end{equation*}
 \begin{equation} \label{e:Estimate of R}
 \le C\big(\lVert\nabla\theta\rVert_{C_tH_x^{s-1}}^2+\lVert\omega\rVert_{C_tH_x^{s-1}}^2+\lVert J\rVert_{C_tH_x^{s-1}}^2 \big)^{1/2}.    
 \end{equation}
 Therefore eventually we have the following estimate for $\zeta^{\pm}$:
 $$  \lVert\zeta^{\pm}\rVert_{C_tH_x^{s-1}}$$
 \begin{equation} \label{e:ZetaEs}
  \le C[1+C\exp(CTU)]^{\lfloor s\rfloor+1}(\lVert \zeta_0^{\pm}\rVert_{H^{s-1}}+CT[1+C\exp(CTU)]^{\lfloor s\rfloor+1}\lVert R^{\pm}\rVert_{C_tH_x^{s-1}})    .
 \end{equation}
 Moreover by lemma 2.2, we have
 \begin{equation}\label{e:Estheta}
  \lVert\nabla\theta\rVert_{C_tH_x^{s-1}}\le C[1+C\exp(CTU)]^{\lfloor s\rfloor+1}\lVert \nabla\theta_0\rVert_{H^{s-1}} .  
 \end{equation}
 \item \textbf{Step 3.} Define now the function space $\mathcal X$
\begin{equation*}
 \mathcal X^s=C([0,T];H_\sigma^s\times H^s\times H_\sigma^s)   
\end{equation*}
 and its $U$-ball $\mathscr B_U$:
 \begin{equation}
\mathscr B_U:=\big\{(u,\theta,B)\in C([0,T];H_\sigma^s\times H^s\times H_\sigma^s):\lVert (u,\theta,B)\rVert_{\mathcal X^s}\le U\big\}     
 \end{equation}
 where $U$ is a positive number, $T$ is sufficiently small positive number to be chosen later, semi-norm $\lVert\cdot\rVert_{\mathcal X}$ is define as:
 \begin{equation}
\lVert (u,\theta,B)\rVert_{\mathcal X^s}=\big(\lVert\nabla u\rVert_{C_tH_x^{s-1}}^2+\lVert\nabla \theta\rVert_{C_tH_x^{s-1}}^2+\lVert\nabla B\rVert_{C_tH_x^{s-1}}^2\big)^{1/2} .    
 \end{equation}
 Now we consider the following iterative scheme: For given 
\begin{equation*}
(u^n,\theta^n,B^n)\in\mathcal X    
\end{equation*}
 we set:
\begin{equation}\label{Picard System}
\left\{
\begin{aligned}
&\dot X^n=u^n(t,X^n),\quad A^n=(X^n)^{-1}\\
&G^n=\theta^n e_d+B^n\cdot\nabla B^n\\
&v^n(t,x)=u_0(x)+\int_0^t\nabla^*X^n(\tau,x)G^n(\tau,X^n(\tau,x))d\tau\\
&u^{n+1}(t,x)=\mathbf P\big[\nabla^*A_t^n(x)v^n(t,A^n(t,x))\big]\\
&\dot X^{n+1}=u^{n+1}(t,X^{n+1}),\quad A^{n+1}=(X^{n+1})^{-1}\\
&\theta^{n+1}=\theta_0\circ A^{n+1}\\
&B^{n+1}=(X^{n+1})_\sharp B_0=(\nabla X^{n+1}B_0)\circ A^{n+1}
\end{aligned}
\right.
\end{equation}
Now by \eqref{Picard System} and \eqref{Div-free Sobolev norm Identity} we have
\begin{equation*}
\lVert(u^{n+1},\theta^{n+1},B^{n+1})\rVert_{\mathcal X}^2=\lVert \omega^{n+1}\rVert_{C_tH_x^{s-1}}^2+\lVert J^{n+1}\rVert_{C_tH_x^{s-1}}^2+\lVert \nabla\theta^{n+1}\rVert_{C_tH_x^{s-1}}^2 .   
\end{equation*}
Thanks to \eqref{e:EnergyEs}\eqref{e:ZetaEs}, denote by $C_{T,s,U}:=[1+C\exp(CTU)]^{\lfloor s\rfloor+1}$ we have
\begin{align*}
&\lVert(u^{n+1},\theta^{n+1},B^{n+1})\rVert_{\mathcal X}^2\\
&\le CC_{T,s,U}(\lVert \zeta_n^{+}\rVert_{H^{s-1}}^2+\lVert \zeta_n^{-}\rVert_{H^{s-1}}^2+\lVert\nabla\theta_0\rVert_{H^{s-1}}^2)\\
&\le CC_{T,s,U}^3[(\lVert \zeta_0^{+}\rVert_{H^{s-1}}+\lVert \zeta_0^{-}\rVert_{H^{s-1}}+CC_{T,s,U}T\lVert R^{\pm}\rVert_{C_tH_x^{s-1}})^2+\lVert\nabla\theta_0\rVert_{H^{s-1}}^2]\\
&\le \underbrace{CC_{T,s,U}^5T^2\lVert R^{\pm}\rVert_{C_tH_x^{s-1}}^2}_{:=(I_1)}+\underbrace{CC_{T,s,U}^3\lVert(u_0,\theta_0,B_0)\rVert_{\mathcal X}^2}_{:=(I_2)}.
\end{align*}
Now we estimate $(I_1), (I_2) $ separately. For $(I_1)$ we notice:
\begin{equation*}
I_1\le  CC_{T,s,U}^5T^2(\lVert (u^n,\theta^n,B^n)\rVert_{\mathcal X}^2).    
\end{equation*}
Choose $T=U^{-1}$ and define $(1+Ce^C)^{\lfloor s\rfloor+1}=C_1$, the above inequality reads as:
\begin{equation*}
I_1\le CC_1^5T^2U^2=CC_1^5 
\end{equation*}
where for $I_2$:
\begin{equation*}
 I_2\le CC_1^3\lVert(u_0,\theta_0,B_0)\rVert_{\mathcal X}^2 .   
\end{equation*}
therefore to ensure the next iteration $(u^{n+1},\theta^{n+1},B^{n+1})$ still lies in $\mathscr B_U$, it suffices to choose $U$ such that
  \begin{equation*}
   U^2\ge I_1+I_2\ge CC_1^3\lVert(u_0,\theta_0,B_0)\rVert_{\mathcal X}^2+CC_1^5 .
  \end{equation*}
  Now existence of solution of above inequality is obvious, essentially we just need to choose $U=\tilde C\lVert (u_0,\theta_0,B_0)\rVert_{\mathcal X}$ for sufficiently large $\tilde C>>1$. With above choice of $U$ we have uniform upper bound for sequence $\{(u^n,\theta^n,B^n)\}_{n=1}^\infty\subset\mathcal X$ as follows:
  \begin{equation}\label{UnifB}
   \sup_{n\in\mathbb N}  \lVert(u^n,\theta^n,B^n)\rVert_{\mathcal X}\le U
  \end{equation}
 \item \textbf{Step 4.}  We now aim to show that $(u^n,\theta^n,B^n)$ is a Cauchy sequence in product space $C_t(L_x^2\times L_x^2\times L_x^2)$. First we have the following estimates for $\theta^n$  by Lemma 2.3:
$$\lVert \theta^{m+1}(t)-\theta^{n+1}(t)\rVert_2$$
$$=\lVert \theta_{0}\circ A^{m+1}(t)-\theta_0\circ A^{n+1}(t)\rVert_2$$
$$=\bigg\lVert\int_0^1\frac{d}{dc} \theta_0(cA^{m+1}(t)+(1-c)A^{n+1}(t))dc \bigg\rVert_2$$
$$\le \lVert\nabla\theta_0\rVert_\infty\lVert A^{m+1}(t)-A^{n+1}(t)\rVert_2$$
\begin{equation}\label{CaucTheta}
\le \lVert\nabla\theta_0\rVert_{\infty}\int_0^t\lVert u^{m+1}-u^{n+1}\rVert_2(\tau) d\tau.    
\end{equation}
Now, to estimate $u^n$, $B^n$ , we indeed consider the sequence $\zeta^{+,n}$ and $\zeta^{-.n}$ defined as:
\begin{equation}\label{zeta ite}
\zeta^{\pm,n+1}=(Y_t)_\sharp \big[\zeta_0^{\pm}+\int_0^t (Z_\tau\circ A_\tau^n)_\sharp \big(\nabla\theta_\tau^n\times e_d\pm \nabla\times M_\tau^n\big)d\tau\big].   
\end{equation}
Here $Y$ is the flow of $-B_0$, $Z$ is the inverse flow of $Y$ and thus $Y$, $Z$ are independent of $n$. Now $u^n$, $B^n$ are recovered from $\zeta^{\pm,n}$ via
\begin{equation*}
u^n=\frac{1}{2}{\BioS}(\zeta^{+,n}+\zeta^{-,n}),\quad B^n=\frac{1}{2}{\BioS}(\zeta^{+,n}-\zeta^{-,n}).    
\end{equation*}
Recalling \eqref{e:Diffrence Weber Sobolev Estimate}, we notice that
\begin{equation*}
\mathcal W(u\circ \ell,\ell)= {\BioS}(\ell_\sharp(\nabla\times\omega))    
\end{equation*}
hence estimate \eqref{e:Difference Weber L2 Estimate} reads:
\begin{equation}\label{zetadiff} 
   \lVert (\ell_1)_\sharp\omega_1-(\ell_2)_\sharp\omega_2\rVert_{H^{-1}}\le C\lVert\omega_1\rVert_{H^{s-1}}\lVert \ell_1-\ell_2\rVert_2+C\lVert\nabla \ell_2\rVert_\infty\lVert \omega_1-\omega_2\rVert_{H^{-1}}.
    \end{equation} 
Now we aim to show that $\zeta^{\pm,n}$ are Cauchy in $H^{-1}$. One first notice since $Y_t$ is $n$-independent, by Lemma 2.1 and measure-preserving property of $Y,Z$ we have the following estimate
$$\lVert \zeta^{\pm,n}-\zeta^{\pm,m}\rVert_{H^{-1}}=\lVert (Y_t)_\sharp[(Z_t)_\sharp(\zeta^{\pm,n}-\zeta^{\pm,m})]\rVert_{H^{-1}}$$
$$=\lVert{\BioS}((Y_t)_\sharp[(Z_t)_\sharp(\zeta^{\pm,n}-\zeta^{\pm,m})])\rVert_{2}$$
$$=\lVert\nabla^*Z_t({\BioS}[(Z_t)_\sharp(\zeta^{\pm,n}-\zeta^{\pm,m})])\circ Z_t\rVert_{2}$$
$$\le \lVert\nabla^*Z_t\rVert_\infty \lVert({\BioS}[(Z_t)_\sharp(\zeta^{\pm,n}-\zeta^{\pm,m})])\circ Z_t\rVert_{2}$$
$$\le \exp(Ct\lVert \nabla B_0\rVert_{H^{s-1}})\lVert [(Z_t)_\sharp(\zeta^{\pm,n}-\zeta^{\pm,m})]\rVert_{H^{-1}}.$$
Then by \eqref{zeta ite} \eqref{zetadiff} we have
\begin{equation*}
\lVert \zeta^{\pm,n+1}-\zeta^{\pm,m+1}\rVert_{H^{-1}}  
\end{equation*}
$$\le \exp(Ct\lVert \nabla B_0\rVert_{H^{s-1}})\lVert [(Z_t)_\sharp(\zeta^{\pm,n+1}-\zeta^{\pm,m+1})]\rVert_{H^{-1}}$$
$$\le\exp(Ct\lVert \nabla B_0\rVert_{H^{s-1}})\int_0^t\lVert (Z_\tau\circ A_\tau^n)_\sharp R_n-(Z_\tau\circ A_\tau^m)_\sharp R_m\rVert_{H^{-1}}d\tau.$$
Denoting by ${\tilde A}^n:=Z\circ A^n$. Now apply \eqref{zetadiff}, choose $\ell_1={\tilde A}^n$, $\ell_2={\tilde A}^m$ and $\omega_1=R_m^\pm$, $\omega_2=R_n^\pm$ we obtain
\begin{equation*}
\lVert \zeta^{\pm,n+1}-\zeta^{\pm,m+1}\rVert_{H^{-1}}   
\end{equation*}
$$\le\exp(Ct\lVert\nabla B_0\rVert_{H^{s-1}})\int_0^t\lVert ({\tilde A}^n)_\sharp R_n^{\pm}-({\tilde A}^m)_\sharp R_m^{\pm}d\tau\rVert_{H^{-1}}d\tau$$
$$\le C\exp(Ct\lVert\nabla B_0\rVert_{H^{s-1}})\bigg[\int_0^t \lVert R_n^{\pm}\rVert_{H^{s-1}}\lVert{\tilde A}^n-{\tilde A}^m\rVert_2 d\tau+\int_0^t \lVert\nabla {\tilde A}^n\rVert_\infty \lVert R_n^{\pm}-R_m^{\pm}\rVert_{H^{-1}}d\tau\bigg ].$$
Now $R_n^\pm$ is bounded in $H^{s-1}$ by $CU^2$. Moreover note that
\begin{equation*}
{\tilde A}=Z\circ A=(X\circ Y)^{-1}    
\end{equation*}
where differentiate $X\circ Y$ in time we have:
\begin{equation*}
 \frac{d}{dt}(X\circ Y)=u\circ (X\circ Y)-(\nabla X B_0)\circ Y=(u-B)\circ(X\circ Y).   
\end{equation*}
Hence, we have $\lVert\nabla\tilde A\rVert_\infty\le \exp(CUT)$. For ${\tilde A}^n-{\tilde A}^m$, Apply lemma 2.3 we obtain
\begin{equation*}
\lVert{\tilde A}^n-{\tilde A}^m\rVert_2\le \exp(CUt)\int_0^t\lVert u^n-u^m\rVert_2+\lVert B^n-B^m\rVert_2 d\tau .   
\end{equation*}
\item\textbf{Step 5.}
We now prove the following estimate for $\lVert R_n^{\pm}-R_m^{\pm}\rVert_{H^{-1}}$ :
\begin{equation}
 \lVert R_n^{\pm}-R_m^{\pm}\rVert_{H^{-1}}\le C\big(\lVert u^n-u^m\rVert_2+\int_0^t\lVert u^m-u^n\rVert_2d\tau+\lVert B^n-B^m\rVert_2\big).   
\end{equation}
Recall $R_n$ is defined as 
\begin{equation*}
R_n^{\pm}=\nabla\times\big(x_d\nabla\theta^n\pm 2\nabla^*u^nB^n\big) .   
\end{equation*}
We use duality to compute $\lVert R_n-R_m\rVert_{H^{-1}}$. Note that
\begin{equation*}
\lVert R_n-R_m\rVert_{H^{-1}}=\sup_{\lVert\psi\rVert_{H^1}\le 1}\lvert\langle \psi, R_n-R_m\rangle_{2}\rvert   
\end{equation*}
\begin{equation*}
\le\underbrace{\sup_{\lVert\psi\rVert_{H^1}\le 1}\lvert\langle \nabla\times\psi, x_d(\nabla\theta^n-\nabla\theta^m)\rangle_{2}\rvert}_{:=I_1}+ 2\underbrace{\sup_{\lVert\psi\rVert_{H^1}\le 1}\lvert\langle \nabla\times \psi, \nabla^*u^nB^n-\nabla^*u^mB^m\rangle_{2}\rvert}_{:=I_2}.    
\end{equation*}
For $I_1$, integrating by part we have the following bound:
\begin{equation*}
\lvert\langle \nabla\times \psi, x_d(\nabla\theta^n-\nabla\theta^m)\rangle_{2}\rvert=\lvert\langle\theta^n-\theta^m,\nabla\times\psi \cdot e_d\rangle_2\rvert\le \lVert\psi\rVert_{H^1}\lVert\nabla\theta_0\rVert_\infty \lVert A^n-A^m\rVert_2.    
\end{equation*}
Therefore we conclude:
\begin{equation*}
I_1\le \sup_{\lVert\psi\rVert_{H^{1}\le 1}} C\lVert\psi\rVert_{H^1}\lVert\nabla\theta_0\rVert_\infty\lVert A^n-A^m\rVert_2\le C\lVert\nabla\theta_0\rVert_{H^{s-1}}\lVert A^n-A^m\rVert_2.   
\end{equation*}
Now we estimate $I_2$. Notice:
\begin{equation*}
 I_2=2\sup_{\lVert\psi\rVert_{H^{1}}\le 1}\lvert\langle \nabla\times\psi, \nabla^*(u^n-u^m)B^n+\nabla^*u^m(B^n-B^m)\rangle_{2}\rvert .  
\end{equation*}
For the first term, integrate by part and use divergence free property of $\nabla\times\psi$:
\begin{equation*}
\lvert\langle \nabla\times\psi, \nabla^*(u^n-u^m)B^n\rangle_2\rvert=\lvert\langle \nabla\times\psi\cdot\nabla B^n, u^n-u^m\rangle_2\rvert \le C\lVert \nabla B^n\rVert_{\infty}\lVert\nabla\times\psi\rVert_2\lVert u^n-u^m\rVert_2.  
\end{equation*}
For the second term we simply estimate it as
\begin{equation*}
\lvert\langle \nabla\times\psi, \nabla^*u^m(B^n-B^m)\rangle_2\rvert\le C\lVert\nabla\times\psi\rVert_2\lVert\nabla u^m\rVert_\infty \lVert B^n-B^m\rVert_2.    
\end{equation*}
Summing up the above estimates we conclude:
$$I_2\le \sup_{\lVert\psi\rVert_{H^{1}}=1}C(\lVert \nabla B^n\rVert_{\infty}\lVert\nabla\times\psi\rVert_2\lVert u^n-u^m\rVert_2+\lVert\nabla\times\psi\rVert_2\lVert\nabla u^m\rVert_\infty \lVert B^n-B^m\rVert_2 )$$
$$\le C(\lVert\nabla u^m\rVert_{H^{s-1}}+\lVert\nabla B^n\rVert_{H^{s-1}})(\lVert u^n-u^m\rVert_2+\lVert B^n-B^m\rVert_2).$$
Therefore, putting above estimates together we conclude:
\begin{equation}\label{R Estimate}
\lVert R_n-R_m\rVert_{H^{-1}}\le C(\lVert u^n-u^m\rVert_2+\lVert B^n-B^m\rVert_2+\int_0^t\lVert u^n-u^m\rVert_2d\tau). 
\end{equation}
\item \textbf{Step 6.}  Thanks to \eqref{R Estimate}, estimates of $\lVert \zeta^{\pm,n+1}-\zeta^{\pm,m+1}\rVert_{H^{-1}}$ goes as follows:
$$\lVert \zeta^{\pm,n+1}-\zeta^{\pm,m+1}\rVert_{H^{-1}}\le CU^2\exp(CUT)\int_0^t \lVert{\tilde A}^n-{\tilde A}^m\rVert_2+\lVert R_n^{\pm}-R_m^{\pm}\rVert_{H^{1}} d\tau$$
\begin{equation}\label{e:zetaest}
\le CU^2\exp(CUT)\bigg[C\int_0^t\lVert u^n-u^m\rVert_2+\lVert B^n-B^m\rVert_2d\tau+\int_0^t\int_0^\tau \lVert u^n-u^m\rVert_2+\lVert B^n-B^m\rVert_2 dr.     
\end{equation}
Now with \eqref{e:zetaest} we are able to estimate:
$$\lVert u^{n+1}-u^{m+1}\rVert_2^2+\lVert B^{n+1}-B^{m+1}\rVert_2^2$$
$$=\frac{1}{2}\big(\lVert (u^{n+1}-u^{m+1})+(B^{n+1}-B^{m+1})\rVert_2^2+\lVert (u^{n+1}-u^{m+1})-(B^{n+1}-B^{m+1})\rVert_2^2\big)$$
$$=\frac{1}{2}\big(\lVert \zeta^{+,n+1}-\zeta^{+,m+1}\rVert_{H^{-1}}^2+\lVert\zeta^{-,n+1}-\zeta^{-,m+1}\rVert_{H^{-1}}^2\big)$$
$$\le \frac{1}{2} CU^4\bigg(\int_0^t \lVert u^{n}-u^{m}\rVert_2^2+\lVert B^{n}-B^{m}\rVert_2^2d\tau+\int_0^t\int_0^\tau \lVert u^{n}-u^{m}\rVert_2^2+\lVert B^{n}-B^{m}\rVert_2^2dr\bigg).$$ 
\textbf{Step 7.} Now set 
$$F_{n+1,m+1}(t):=\lVert u^{n+1}-u^{m+1}\rVert_2^2+\lVert B^{n+1}-B^{m+1}\rVert_2^2,\quad \tilde F:=\limsup_{n,m\to\infty} F_{n,m}.$$ 
By the above estimate we have:
\begin{equation}\label{Gron}
\tilde F(t)\le \int_0^t \tilde F(\tau)d\tau+\int_0^t\int_0^\tau \tilde F(s)ds. 
\end{equation}
Apply the Gronwall's inequality to \eqref{Gron}, together with $\tilde F(0)=0$ we conclude that $\tilde F=0$. Therefore, both $u^n$ and $B^n$ are Cauchy sequence in $L^2$, by \eqref{CaucTheta} we conclude $\theta^n$ is also Cauchy in $L^2$. Now thanks to \eqref{UnifB}, by interpolation we therefore obtain for any $\varepsilon>0$:
\begin{equation*}
\limsup_{m,n\to\infty}\lVert (u^n,\theta^n,B^n)-(u^m,\theta^m,B^m)\rVert_{\mathcal X^{s-\varepsilon}}=0 .   
\end{equation*}
Therefore, there exists a unique limit $(u^*,\theta^*,B^*)$ solving \eqref{Picard System}, which also lies in $\mathscr B_U$ by compactness. Hence, the proof is completed. \qed
\end{itemize}
\paragraph{Remark 4.2} We remark that the proof of local well-posedness of Camassa Holm/ GSQG type equation in Sobolev class is also available, following the same strategy. In which case we simply consider Fourier multiplier $\bT$ and potential $\sP$ satisfies the following:
\begin{equation*}
\bT: H^{s}\to H^{s+\beta}, \quad -2\le \beta\le 1;\quad \sP=0    
\end{equation*}
Here the borderline case  $\bT=(-\Delta)^{-1/2}$ meets the critical SQG exponent $1/2$, and Lagrangian approach may fail for supercritical exponent. We omit the details here and refer the readers to \cite{Pa26}. 

\paragraph{Remark 4.3} It's not applicable to work with very general $\sP$ in \eqref{e:Stochastic Lagrangian System}, even for $\bT=\mathbf{Id}$, as one can observe in standard local well-posedness proof of MHD, when one apply Kato-Ponce commutator estimate, the cancellation of bad term crucially relies on the magnetic convection term:
\begin{equation*}
 \langle\Lambda^s(B\cdot\nabla B),\Lambda^s u\rangle_2+\langle \Lambda^s(B\cdot\nabla u),\Lambda^s B\rangle_2   
\end{equation*}
\begin{equation*}
\le \langle[\Lambda^s,B\cdot\nabla]B,\Lambda^s u\rangle+\langle[\Lambda^s,B\cdot\nabla]u,\Lambda^s B\rangle_2-(\langle B\cdot\nabla\Lambda^s \nabla B,\Lambda^s u\rangle_2+\langle B\cdot\nabla\Lambda^s \nabla u,\Lambda^s B\rangle_2)    
\end{equation*}
\begin{equation*}
 \le C(\lVert \nabla u\rVert_{H^{s}}^2+\lVert \nabla B\rVert_{H^s}^2)-\underbrace{\langle B, \nabla(\Lambda^su\cdot\Lambda^sB)\rangle_2}_{=0}.   
\end{equation*}
Assume $\sP(X)=P_1(X_\sharp\theta_0)+P_2(X_\sharp B_0)$ with $P_1: H^s\to\mathbb R$, $P_2:H_\sigma^s\to \mathbb R$. Then general choice of $P_2$ will deteriorate above cancellation, which will also fail the fixed point argument. However, if $\delta P_2/\delta w(B)$ is smoother than $B$, or
\begin{equation*}
\frac{\delta P_2}{\delta w}(B)-B\in H_\sigma^{s}\quad\text{provided } B\in H_\sigma^s.    
\end{equation*}
Then the above theorem still holds true.

\bibliographystyle{siam}  
\bibliography{refs}

\end{document}